\documentclass[12pt]{article}

\usepackage{amsbsy,amsmath,amsthm,amssymb}

\usepackage[sort,longnamesfirst]{natbib}
\newcommand{\pcite}[1]{\citeauthor{#1}'s \citeyearpar{#1}}
\usepackage{geometry}
\usepackage[dvips]{changebar}

\theoremstyle{remark}
\newtheorem{remark}{Remark}
\newtheorem{example}{Example}
\newtheorem{lemma}{Lemma}
\newtheorem{corollary}{Corollary}
\newtheorem*{theorem}{Theorem}
\newtheorem*{proposition}{Proposition}

\newcommand{\gbar}{\bar g}
\newcommand{\sX}{{\mathsf X}}

\begin{document}

\title{Fixed-Width Output Analysis for Markov Chain Monte Carlo}
\author{Galin L. Jones \\ School of Statistics \\ University of
  Minnesota \\ {\tt galin@stat.umn.edu} \and Murali Haran \\Department
  of Statistics \\ The Pennsylvania State University \\ {\tt
    mharan@stat.psu.edu} \and Brian S. Caffo \\ Department of
  Biostatistics\\ Johns Hopkins University \\ {\tt
    bcaffo@jhsph.edu} \and Ronald Neath \\ School of Statistics \\
  University of Minnesota \\ {\tt rneath@stat.umn.edu}}
\date{February 9, 2005 (Revised October 5, 2005)}\maketitle

\begin{abstract}
  Markov chain Monte Carlo is a method of producing a correlated
  sample in order to estimate features of a target distribution via
  ergodic averages.  A fundamental question is when should sampling
  stop?  That is, when are the ergodic averages good estimates of the
  desired quantities?  We consider a method that stops the simulation
  when the width of a confidence interval based on an ergodic average
  is less than a user-specified value.  Hence calculating a Monte
  Carlo standard error is a critical step in assessing the simulation
  output.  We consider the regenerative simulation and batch means
  methods of estimating the variance of the asymptotic normal
  distribution.  We give sufficient conditions for the strong
  consistency of both methods and investigate their finite sample
  properties in a variety of examples.
\end{abstract}

\section{Introduction}
\label{sec:intro}
Suppose our goal is to calculate $E_{\pi} g := \int_{\mathsf{X}} g (x)
\pi(dx)$ with $\pi$ a probability distribution having support $\sX$
and $g$ a real-valued, $\pi$-integrable function.  Also, suppose $\pi$
is such that Markov chain Monte Carlo (MCMC) is the only viable method
for estimating $E_{\pi} g$.

Let $X = \{X_0, X_1, X_2, \dots\}$ be a time-homogeneous, aperiodic,
$\pi$-irreducible, positive Harris recurrent Markov chain with state
space $(\sX, \cal{B} (\sX))$ and invariant distribution $\pi$.  (See
\citet{meyn:twee:1993} for definitions.)  In this case, we say that
$X$ is Harris ergodic and the Ergodic Theorem implies that, with
probability 1,
\begin{equation}
\label{eq:erg_avg}
\bar{g}_{n} := \frac{1}{n} \sum_{i=0}^{n-1} g
(X_i) \rightarrow E_{\pi} g \quad\text{as $n \rightarrow \infty$.}  
\end{equation}
Given an MCMC algorithm that simulates $X$ it is conceptually easy to
generate large amounts of data and use $\bar{g}_{n}$ to obtain an
arbitrarily precise estimate of $E_{\pi}g$.

There are several methods for deciding when $n$ is sufficiently large;
i.e., when to terminate the simulation.  The simplest is to terminate
the computation whenever patience runs out.  This approach is
unsatisfactory since the user would not have any idea about the
accuracy of $\bar{g}_{n}$.  Alternatively, with several preliminary
(and necessarily short) runs the user might be able to make an
informed guess about the variability in $\bar{g}_{n}$ and hence make
an a priori choice of $n$.  Another method would be to monitor the
sequence of $\bar{g}_{n}$ until it appears to have stabilized.  None
of these methods are automated and hence are inefficient uses of user
time and Monte Carlo resources.  Moreover, they provide only a point
estimate of $E_{\pi} g$ without additional work.

Convergence diagnostics are also sometimes used to terminate the
simulation \citep{cowl:carl:1996}.  Some convergence diagnostics are
available in software, e.g. the R package {\tt boa}, and hence may be
considered automated.  However, none of the diagnostics of which we
are aware explicitly address how well $\bar{g}_{n}$ estimates $E_{\pi}
g$; this is discussed again in subsection~\ref{sec:compare}.

An alternative is to calculate a Monte Carlo standard error and use it
to terminate the simulation when the width of a confidence interval
falls below a specified value.  Under regularity conditions (see
Section~\ref{sec:mc_theory}) the Markov chain $X$ and function $g$
will admit a central limit theorem (CLT); that is,
\begin{equation}
\label{eq:clt}
\sqrt{n} (\bar{g}_{n} - \text{E}_{\pi} g) \stackrel{d}{\rightarrow}
\text{N} (0, \sigma^{2}_{g})
\end{equation}
as $n \rightarrow \infty$ where $\sigma^{2}_{g} := \text{var}_{\pi} \{
g(X_{0})\} + 2 \sum_{i=1}^{\infty} \text{cov}_{\pi} \{ g(X_{0}),
g(X_{i})\}$.  Given an estimate of $\sigma_{g}^{2}$, say
$\hat{\sigma}_{n}^{2}$, we can form a confidence interval for
$\text{E}_{\pi} g$.  If this interval is too large then the value of
$n$ is increased and simulation continues until the interval is
sufficiently small; this is a common way of choosing $n$ \citep[e.g.,
see][]{fish:1996,geye:1992,jone:hobe:2001}.  Notice that the final
Monte Carlo sample size is random.  We study sequential fixed-width
methods which formalize this approach.  In particular, the simulation
terminates the first time
\begin{equation}
\label{eq:hw}
t_{*}\,  \frac{\hat{\sigma}_{n}}{\sqrt{n}} + p(n) \le \epsilon
\end{equation}
where $t_{*}$ is an appropriate quantile, $p(n) \ge 0$ on
$\mathbb{Z}_{+}$ and $\epsilon >0$ is the desired half-width.  The
role of $p$ is to ensure that the simulation is not terminated
prematurely due to a poor estimate of $\sigma_{g}^{2}$.  One
possibility is to fix $n^{*} > 0$ and take $p(n) = \epsilon I(n\le
n^{*})$ where $I$ is the usual indicator function.

Sequential statistical procedures have a long history; see
\citet{lai:2001} for an overview and commentary.  Moreover, classical
approaches to sequential fixed-width confidence intervals such as
those found in \citet{chow:robb:1965}, \citet{liu:1997} and
\citet{nada:1969} are known to work well.  However, the classical
procedures are not relevant to the current work since they assume the
observations are random samples.

In a simulation context, procedures based on \eqref{eq:hw} were
studied most notably by \citet{glyn:whit:1992} who established that
these procedures are \textit{asymptotically valid} in that if our goal
is to have a $100(1-\delta)\%$ confidence interval with width
$2\epsilon$ then
\begin{equation}
\Pr (\text{E}_{\pi} g \, \in \, \text{Int}[T(\epsilon)]) \rightarrow 1
- \delta \hspace*{3mm} \text{ as } \, \epsilon \rightarrow 0
\end{equation}
where $T(\epsilon)$ is the first time that \eqref{eq:hw} is satisfied
and $\text{Int}[T(\epsilon)]$ is the interval at this time.  Glynn and
Whitt's conditions for asymptotic validity are substantial: (i) A
functional central limit theorem (FCLT) holds; (ii)
$\hat{\sigma}^{2}_{n} \rightarrow \sigma_{g}^{2}$ with probability 1
as $n \rightarrow \infty$; and (iii) $p(n)=o(n^{-1/2})$.  Markov
chains frequently enjoy an FCLT under the same conditions that ensure
a CLT.  However, in the context of MCMC, little work has been done on
establishing conditions for (ii) to hold.  Thus one of our goals is to
give conditions under which some common methods provide strongly
consistent estimators of $\sigma_{g}^{2}$.  Specifically, our
conditions require the sampler to be either uniformly or geometrically
ergodic.  The MCMC community has expended considerable effort in
establishing such mixing conditions for a variety of samplers; see
\citet{jone:hobe:2001} and \citet{robe:rose:1998b,robe:rose:2004} for
some references and discussion.

We consider two methods for estimating the variance of the asymptotic
normal distribution, regenerative simulation (RS) and non-overlapping
batch means (BM).  Both have strengths and weaknesses; essentially, BM
is easier to implement but RS is on a stronger theoretical footing.
For example, when used with fixed number of batches BM \textit{cannot}
be even weakly consistent for $\sigma_{g}^{2}$.  We give conditions
for the consistency of RS and show that BM can provide a consistent
estimation procedure by allowing the batch sizes to increase (in a
specific way) as $n$ increases.  In this case it is denoted CBM to
distinguish it from the standard fixed-batch size version which we
denote BM.  This was addressed by \citet{dame:1994} but, while the
approach is similar, our regularity conditions on $X$ are weaker.
Also, the regularity conditions required to obtain strong consistency
of the batch means estimator are slightly stronger than those required
by RS.  Finally, it is important to note that RS and CBM do not
require that $X$ be stationary; hence burn-in is not required.

The justification of fixed-width methods is entirely asymptotic so it
is not clear how the finite sample properties of BM, CBM, and RS
compare in typical MCMC settings.  For this reason, we conduct a
simulation study in the context of two benchmark examples and two
realistic examples, one of which is a complicated frequentist problem
and one which involves a high-dimensional posterior.  Roughly
speaking, we find that BM performs poorly while RS and CBM are
comparable.  

The rest of this article is organized as
follows. Section~\ref{sec:mc_theory} fixes notation and contains a
brief discussion of some relevant Markov chain theory.  In
Section~\ref{sec:oa} we consider RS and CBM.  Then
Section~\ref{sec:examples} contains the examples. 

\section{Basic Markov Chain Theory}
\label{sec:mc_theory}
For $n \in \mathbb{N} := \{1,2,3,\ldots \}$ let $P^n(x,dy)$ be the
$n$-step Markov transition kernel; that is, for $x \in \mathsf{X}$ and
$A \in \cal{B} (\sX)$, $P^n(x,A) = \Pr\left(X_n \in A|X_0 = x\right)$.
A Harris ergodic Markov chain $X$ enjoys a strong form of convergence.
Specifically, if $\lambda(\cdot)$ is a probability measure on
$\cal{B}(\sX)$ then
\begin{equation}
  \label{eq:bas_con}
  \| P^n(\lambda,\cdot) - \pi(\cdot)\| \; \downarrow \; 0
  \quad\text{as $n \rightarrow \infty,$}
\end{equation}
where $P^n(\lambda,A) := \int_{\sX} P^{n} (x, A) \lambda(dx)$ and
$\|\cdot\|$ is the total variation norm.  Suppose there exists an
extended real-valued function $M(x)$ and a nonnegative decreasing
function $\kappa(n)$ on $\mathbb{Z}_{+}$ such that
\begin{equation}
\label{eq:tvbd}
\| P^{n} (x, \cdot) - \pi(\cdot) \| \le M(x) \kappa(n) \; .
\end{equation}
When $\kappa(n) = t^{n}$ for some $t < 1$ say $X$ is
\textit{geometrically ergodic} if $M$ is unbounded and
\textit{uniformly ergodic} if $M$ is bounded.  \textit{Polynomial
  ergodicity of order m} where $m \ge 0$ means $M$ may be unbounded
and $\kappa(n)=n^{-m}$.  

Also, $P$ satisfies \textit{detailed balance} with respect to $\pi$ if
\begin{equation}
\label{eq:dbc}
\pi(dx) P(x,dy) = \pi(dy) P(y, dx) \hspace{5mm} \text{ for all } \,
x,y \in \sX \; .
\end{equation}
Note that Metropolis-Hastings samplers satisfy \eqref{eq:dbc} by
construction but many Gibbs samplers do not.  We are now in position
to give conditions for the existence of a CLT.

\begin{theorem}
\label{thm:clt}
Let $X$ be a Harris ergodic Markov chain on $\mathsf{X}$ with
invariant distribution $\pi$ and suppose $g : \mathsf{X} \rightarrow
\mathbb{R}$ is a Borel function. Assume one of the following
conditions: \vspace{-2mm}
\begin{enumerate}
\item $X$ is polynomially ergodic of order $m > 1$, $\text{E}_{\pi} M
  < \infty$ and there exists $B< \infty$ such that $|g(x)| < B$ almost
  surely; \vspace{-1mm}
\item $X$ is polynomially ergodic of order $m$, $\text{E}_{\pi} M <
  \infty$ and $\text{E}_{\pi} |g(x)|^{2+\delta} < \infty$ for some
  $\delta > 0$ where $m\delta > 2+\delta$; \vspace{-1mm}
\item $X$ is geometrically ergodic and $\text{E}_{\pi} [g^{2}(x)
  (\log^{+} |g(x)|)] < \infty$;  \vspace{-1mm}
\item $X$ is geometrically ergodic, satisfies \eqref{eq:dbc} and
  $\text{E}_{\pi} g^{2}(x) < \infty$; or \vspace{-1mm}
\item $X$ is uniformly ergodic and $\text{E}_{\pi} g^{2}(x) < \infty$.
\end{enumerate} \vspace{-2mm}
Then, for any initial distribution, as $n \rightarrow \infty$
\[
\sqrt{n} (\bar{g}_{n} - \text{E}_{\pi} g) \stackrel{d}{\rightarrow}
\text{N} (0, \sigma_{g}^{2}) \; .
\]
\end{theorem}

\begin{remark}
  The theorem was proved by \citet{ibra:linn:1971} (condition 5),
  \citet{robe:rose:1997c} (condition 4), \citet{douk:mass:rio:1994}
  (condition 3).  See \citet{jone:2004} for details on conditions 1
  and 2.
\end{remark}

\begin{remark}
  Conditions 3, 4 and 5 of the theorem are also sufficient to
  guarantee the existence of an FCLT; see \citet{douk:mass:rio:1994},
  \citet{robe:rose:1997c} and \citet{bill:1968}, respectively.
\end{remark}

\begin{remark}
  The mixing conditions on the Markov chain $X$ stated in
  Theorem~\ref{thm:clt} are not necessary for the CLT; see, for
  example, \citet{chen:1999}, \citet{meyn:twee:1993} and
  \citet{numm:2002}.  However, the weaker conditions are often
  prohibitively difficult to check in situations where MCMC is
  appropriate.
\end{remark}

\begin{remark}
\label{rm:geo}
There are constructive techniques for verifying the existence of an
appropriate $M$ and $\kappa$ from \eqref{eq:tvbd} \citep[Ch.
15]{meyn:twee:1993}.  For example, one method of establishing
geometric ergodicity requires finding a function $V : \sX \rightarrow
[1,\infty)$ and a small set $C \in{\cal B}(\sX)$ such that
\begin{equation}
\label{eq:drift}
PV(x) \le \lambda V(x) + b I(x \in C) \hspace*{4mm} \forall \; x \in
\sX 
\end{equation} 
where $PV(x) := \int V(y) P(x, dy)$, $0 < \lambda < 1$ and $b <
\infty$.  Substantial effort has been devoted to establishing
convergence rates for MCMC algorithms via \eqref{eq:drift} or related
techniques.  For example, \citet{hobe:geye:1998},
\citet{hobe:jone:pres:rose:2002}, \citet{jone:hobe:2004},
\citet{marc:hobe:2004}, \citet{mira:tier:2002}, \citet{robe:1995},
\citet{robe:pols:1994}, \citet{robe:rose:1999a},
\citet{rose:1995a,rose:1996} and \citet{tier:1994} examined Gibbs
samplers while \citet{chri:moll:waag:2001},
\citet{douc:fort:moul:soul:2004},
\citet{fort:moul:2000,fort:moul:2003}, \citet{geye:1999},
\citet{jarn:hans:2000}, \citet{jarn:robe:2002},
\citet{meyn:twee:1994b}, and \citet{meng:twee:1996} analyzed
Metropolis-Hastings algorithms.
\end{remark}

\subsection{The Split Chain}
\label{sec:split}
An object that is important to the study of both RS and CBM is the
\textit{split chain} $X':= \left \{ (X_0, \delta_0), (X_1, \delta_1),
  (X_2, \delta_2), \dots \right \}$ which has state space $\sX \times
\{0,1\}$.  The construction of $X'$ requires a \textit{minorization
  condition}; i.e., a function $s: \sX \mapsto [0,1]$ for which
$E_{\pi} s > 0$ and a probability measure $ Q$ such that
\begin{equation}
  \label{eq:mc}
  P(x,A) \geq s(x) \, Q(A) \hspace*{5mm} \text{ for all } x \in \sX
  \text{ and } A \in \cal{B} (\sX) \; .
\end{equation}
When $\sX$ is countable it is easy to see that \eqref{eq:mc} holds by
fixing $x_{*} \in \sX$, setting $s(x) = I(x=x_{*})$ and $Q(\cdot) =
P(x_{*}, \cdot)$.  \citet{mykl:tier:yu:1995} and \citet{rose:1995a}
give prescriptions that are often useful for establishing
\eqref{eq:mc} in general spaces.  Note that \eqref{eq:mc} allows us to
write $P(x,dy)$ as a mixture of two distributions,
\[
  P(x,dy) = s(x) \, Q(dy) + \left[ 1-s(x) \right] R(x,dy),
\]
where $R(x,dy) := \left[ 1-s(x) \right]^{-1} \left[P(x,dy) - s(x) \,
  Q(dy) \right]$ is the \textit{residual} distribution (define
$R(x,dy)$ as 0 if $s(x) = 1$).  This mixture gives us a recipe for
simulating $X'$: given $X_i=x$, generate $\delta_i \sim
\text{Bernoulli}(s(x))$.  If $\delta_i=1$, then draw $X_{i+1} \sim
Q(\cdot)$, else draw $X_{i+1} \sim R(x,\cdot)$.

The two chains, $X$ and $X'$ are closely related since $X'$ will
inherit properties such as aperiodicity and positive Harris recurrence
and the sequence $\{X_i : i=0,1,\dots\}$ obtained from $X'$ has the
same transition probabilities as $X$.  Also, $X$ and $X'$ converge to
their respective stationary distributions at exactly the same rate.

If $\delta_i=1$, then time $i+1$ is a \textit{regeneration time} when
$X'$ probabilistically restarts itself.  Specifically, suppose we
start $X'$ with $X_0 \sim Q$.  Then each time that $\delta_i=1$,
$X_{i+1} \sim Q$.  Let $0=\tau_{0} < \tau_{1} < \cdots$ be the
regeneration times.  That is, set $\tau_{r+1} = \min \{ i > \tau_{r}
\, : \, \delta_{i-1}=1\}$.  Also assume that $X'$ is run for $R$
tours; that is, the simulation is stopped the $R$th time that a
$\delta_i=1$.  Let $\tau_{R}$ denote the total length of the
simulation and $N_r$ be the length of the $r$th tour; that is, $N_r =
\tau_r - \tau_{r-1}$.  Define
\begin{equation*}
  S_{r} = \sum_{i=\tau_{r-1}}^{\tau_r-1} g (X_i)
\end{equation*}
for $r=1,\ldots,R$.  The $(N_r,S_{r})$ pairs are iid since each is
based on a different tour.  In the sequel we will make repeated use of
the following lemma which generalizes Theorem 2 of
\citet{hobe:jone:pres:rose:2002}.

\begin{lemma}\label{lemma:rsmoments}
  Let $X$ be a Harris ergodic Markov chain with invariant distribution
  $\pi$.  Assume that \eqref{eq:mc} holds and that $X$ is
  geometrically ergodic.  Let $p \ge 1$ be an integer.
\vspace{-2mm}
\begin{enumerate}
\item If $E_{\pi} |g|^{2^{(p-1)} + \delta} < \infty$ for some $ \delta
  > 0$ then $E_{Q} N_{1}^{p} < \infty$ and $E_{Q} S_{1}^{p} < \infty$.
\vspace{-2mm}
\item If $E_{\pi} |g|^{2^{p} + \delta} < \infty$ for some $\delta > 0$
  then $E_{Q} N_{1}^{p} < \infty$ and $E_{Q} S_{1}^{p+ \delta} <
  \infty$.
\end{enumerate} 
\end{lemma}

\begin{proof}  See Appendix~\ref{app:rsmoments}. \end{proof}

\section{Output Analysis}
\label{sec:oa}

\subsection{Regenerative Simulation}
\label{sec:rs}
Regenerative simulation is based on directly simulating the split
chain.  However, using the mixture approach described above is
problematic since simulation from $R(x, dy)$ is challenging.
\citet{mykl:tier:yu:1995} suggest a method for avoiding this issue.
Suppose \eqref{eq:mc} holds and that the measures $P(x, \cdot)$ and
$Q(\cdot)$ admit densities $k(\cdot | x)$ and $q(\cdot)$,
respectively.  Then the following recipe allows us to simulate $X'$.
Assume $X_{0} \sim q(\cdot)$; this is typically quite easy to do, see
\citet{mykl:tier:yu:1995} for some examples. Also, note that this
means burn-in is irrelevant.  Draw $X_{i+1} \sim
k(\cdot | x)$, that is, draw from the sampler at hand, and get
$\delta_{i}$ by simulating from the distribution of $\delta_{i} |
X_{i}, X_{i+1}$ with
\begin{equation}
\label{eq:reg_pr}
\Pr(\delta_{i} = 1\, | \, X_{i}, X_{i+1}) = \frac{s(X_{i})
  q(X_{i+1})}{k(X_{i+1} \, | \, X_{i})} \; .
\end{equation}

\begin{example}
\label{ex:indep}
In a slight abuse of notation let $\pi$ also denote the density of the
target distribution.  Consider an independence Metropolis-Hastings
sampler with proposal density $\nu$.  This chain works as follows: Let
the current state be $X_{i}=x$.  Draw $y \sim \nu$ and independently
draw $u \sim \text{Uniform} (0,1)$.  If
\[
u < \frac{\pi(y) \nu(x)}{\pi(x) \nu(y)}
\]
then set $X_{i+1}=y$ otherwise set $X_{i+1}=x$.
\citet{mykl:tier:yu:1995} derive \eqref{eq:reg_pr} for this case.  Let
$c > 0$ be a user-specified constant.  Then conditional on an
acceptance, i.e. $X_{i}=x$ and $X_{i+1}=y$
\begin{equation}
\label{eq:indep_rp}
\Pr(\delta_{i} = 1\, | \, X_{i}=x, X_{i+1}=y) = 
\begin{cases}
  c \, \max \left\{ \frac{\nu(x)}{\pi(x)}, \,
    \frac{\nu(y)}{\pi(y)}\right\} & \text{if } \min \left\{
    \frac{\pi(x)}{\nu(x)}, \, \frac{\pi(y)}{\nu(y)}
  \right\} > c \\
  \frac{1}{c} \, \max \left\{ \frac{\pi(x)}{\nu(x)}, \,
    \frac{\pi(y)}{\nu(y)}\right\} & \text{if } \max \left\{
    \frac{\pi(x)}{\nu(x)}, \, \frac{\pi(y)}{\nu(y)} \right\} < c \\
  1 & \text{otherwise} \; .
\end{cases}
\end{equation}
Note that we do not need to know the normalizing constants of $\pi$ or
$\nu$ to calculate \eqref{eq:indep_rp}.
\end{example}
In discrete state spaces regenerations can be easy to identify.  In
particular, a regeneration occurs whenever the chain returns to any
fixed state; for example, when the Metropolis-Hastings chain accepts a
move to the fixed state.  This regeneration scheme is most useful when
the state space is not too large but potentially complicated; see
subsection~\ref{sec:p-value}.  It will not be useful when the state
space is extremely large because returns to the fixed state are too
infrequent.  Further practical advice on implementing and automating
RS is given in \citet{broc:kada:2005}, \citet{gilk:robe:sahu:1998},
\citet{geye:thom:1995}, \citet{hobe:jone:pres:rose:2002},
\citet{hobe:jone:robe:2005} and \citet{jone:hobe:2001}.

Implementation of RS is simple once we can effectively simulate the
split chain. For example, the Ergodic Theorem implies that
\[
\bar{g}_{\tau_{R}} = \frac{1}{\tau_{R}} \sum_{j=0}^{\tau_{R} - 1}
g(X_{j}) \rightarrow \text{E}_{\pi} g
\]
with probability 1 as $R \rightarrow \infty$ and hence estimating
$\text{E}_{\pi} g$ is routine.

We now turn our attention to calculating a Monte Carlo standard error
for $\bar{g}_{\tau_{R}}$.  Let $E_{Q}$ denote the expectation for the
split chain started with $X_{0} \sim Q(\cdot)$.  Also, let ${\bar N}$
be the average tour length; that is, ${\bar N}=R^{-1} \sum_{r=1}^R
N_{r}$.  Since the $(N_r,S_{r})$ pairs are iid the strong law implies
with probability 1, $\bar{N} \rightarrow E_{Q} N_{1}$ which is finite
by positive recurrence.  If $E_Q N_1^2 < \infty$ and $E_Q S_{1}^2 <
\infty$ it follows that a CLT holds; i.e., as $R \rightarrow \infty$
\begin{equation}
\label{eq:rs_clt}
\sqrt{R} (\bar{g}_{\tau_{R}} - \text{E}_{\pi} g)
\stackrel{d}{\rightarrow} \text{N}\,(0,\xi_{g}^{2}) 
\end{equation} 
where, as shown in \citet{hobe:jone:pres:rose:2002}, $\xi_{g}^{2} =
E_{Q} (S_{1} - N_{1} E_{\pi} g)^{2} / (E_{Q} N_{1} )^{2}$.  An obvious
estimator of $\xi^{2}_{g}$ is
\begin{equation*}
  \hat{\xi}_{RS}^{2} := \frac{1}{\bar{N}^{2}} \frac{1}{R} \sum_{r=1}^R
  (S_{r} - \gbar_{\tau_{R}} N_{r})^{2}   \; .
\end{equation*}
Now consider
\begin{equation*}
\begin{split}
  \hat{\xi}_{RS}^{2} - \xi^{2}_{g} & = \frac{1}{\bar{N}^{2}}
  \frac{1}{R} \sum_{r=1}^R (S_{r} - \gbar_{\tau_{R}} N_{r})^{2} -
  \frac{E_{Q} (S_{1} - N_{1} E_{\pi} g)^{2}}{(E_{Q} N_{1} )^{2}} \pm
  \frac{E_{Q}
    (S_{1} - N_{1} E_{\pi} g)^{2}}{\bar{N}^{2}} \\
  & = \frac{1}{\bar{N}^{2}} \frac{1}{R} \sum_{r=1}^R \left[(S_{r} -
    \gbar_{\tau_{R}} N_{r})^{2} - E_{Q} (S_{1} - N_{1} E_{\pi} g)^{2}
    \pm (S_{r} - N_{r} E_{\pi}g)^{2} \right]  \\
& \hspace*{25mm} + \left[ E_{Q} (S_{1} - N_{1} E_{\pi} g)^{2} \left(
    \frac{1}{\bar{N}^{2}} - \frac{1}{E_{Q} N_{1}^{2}} \right) \right]
\; .
\end{split}
\end{equation*}
Using this representation and repeated application of the strong law
shows that $\hat{\xi}_{RS}^{2} - \xi^{2}_{g} \rightarrow 0$ with
probability 1 as $R \rightarrow \infty$ \citep[also
see][]{hobe:jone:pres:rose:2002}.  It is typically difficult to check
that $E_Q N_1^2 < \infty$ and $E_Q S_{1}^2 < \infty$.  However, using
Lemma~\ref{lemma:rsmoments} yields the following result.
\begin{proposition}
  Let $X$ be a Harris ergodic Markov chain with invariant distribution
  $\pi$.  Assume that $E_{\pi} |g|^{2+\delta} < \infty$ for some
  $\delta > 0$, \eqref{eq:mc} holds and that $X$ is geometrically
  ergodic.  Then \eqref{eq:rs_clt} holds and $\hat{\xi}_{RS}^{2}
  \rightarrow \xi_{g}^{2}$ w. p. 1 as $R \rightarrow \infty$.
\end{proposition}

Fix $\epsilon > 0$ and let $z$ denote an appropriate standard normal
quantile.  An asymptotically valid fixed-width procedure results by
terminating the simulation the first time
\begin{equation}
\label{eq:rs_fw}
z\,  \frac{\hat{\xi}_{RS}}{\sqrt{R}} + p(R) \le \epsilon \; .
\end{equation}

\subsection{Batch Means}
\label{sec:bm}

In standard batch means the output of the sampler is broken into
batches of equal size that are assumed to be approximately
independent.  (This is not strictly necessary; c.f., the method of
overlapping batch means.)  Suppose the algorithm is run for a total of
$n=ab$ iterations (hence $a=a_{n}$ and $b=b_{n}$ are implicit
functions of $n$) and define 
\begin{equation*}
\bar{Y}_{j} := \frac{1}{b} \sum_{i=(j-1)b}^{jb-1} g (X_{i})
\hspace*{5mm} \text{ for } j=1,\ldots,a \; .
\end{equation*}
The batch means estimate of $\sigma_{g}^{2}$ is
\begin{equation}
\label{eq:bmvar}
\hat{\sigma}_{BM}^{2} = \frac{b}{a-1} \sum_{j=1}^{a} (\bar{Y}_{j} -
\bar{g}_n)^{2}  \; .
\end{equation}
With a fixed number of batches \eqref{eq:bmvar} is not a consistent
estimator of $\sigma_{g}^{2}$ \citep{glyn:igle:1990,glyn:whit:1991}.
On the other hand, if the batch size \textit{and} the number of
batches are allowed to increase as the overall length of the
simulation does it may be possible to obtain consistency.  The first
result in this direction is due to \citet{dame:1994} which we now
describe.  The major assumption made by \citet{dame:1994} is the
existence of a strong invariance principle.  Let $B=\{B(t), t \ge 0\}$
denote a standard Brownian motion.  A strong invariance principle
holds if there exists a nonnegative increasing function $\gamma(n)$ on
the positive integers, a constant $0 < \sigma_{g} < \infty$ and a
sufficiently rich probability space such that
\begin{equation}
\label{eq:sip}
\left| \sum_{i=1}^{n} g(X_{i}) - n \text{E}_{\pi} g - \sigma_{g} B(n)
\right| = O(\gamma(n)) \hspace*{5mm} \text{w.p. 1 as } \; n
\rightarrow \infty \;  
\end{equation}
where the w.p. 1 in \eqref{eq:sip} means for almost all sample paths.
In particular, \citet{dame:1994} assumed \eqref{eq:sip} held with
$\gamma(n) = n^{1/2 - \alpha}$ where $ 0 < \alpha \le 1/2$.  However,
it would seem a daunting task to directly check this condition in any
given application.  In an attempt to somewhat alleviate this
difficulty we have the following lemma.

\begin{lemma}  \label{lemma:sip}
  Let $g : \sX \rightarrow \mathbb{R}$ be a Borel function and let
  $X$ be a Harris ergodic Markov chain with invariant distribution
  $\pi$.\vspace{-2mm}
\begin{enumerate}
\item If $ X$ is uniformly ergodic and $\text{E}_{\pi} |g|^{2 +
    \delta} < \infty$ for some $\delta > 0$ then \eqref{eq:sip} holds
  with $\gamma(n) = n^{1/2 - \alpha}$ where $\alpha < \delta/(24+12
    \delta)$. \vspace{-2mm}
\item If $X$ is geometrically ergodic, \eqref{eq:mc} holds and
  $\text{E}_{\pi} |g|^{4 + \delta} < \infty$ for some $\delta > 0$
  then \eqref{eq:sip} holds with $\gamma(n) = n^{\alpha} \log n$ where
  $\alpha = 1/(2 + \delta)$.\vspace{-2mm}
\end{enumerate}
\end{lemma}

\begin{proof}
  The first part of the lemma is an immediate consequence of Theorem
  4.1 of \citet{phil:stou:1975} and the fact that uniformly ergodic
  Markov chains enjoy exponentially fast uniform mixing. The second
  part follows from our Lemma~\ref{lemma:rsmoments} and Theorem 2.1 in
  \citet{csak:csor:1995}.
\end{proof}
  
Using part 1 of Lemma~\ref{lemma:sip} we can state Damerdji's result
as follows. 

\begin{proposition} \label{prop:ue_sc_bm} \citep{dame:1994} 
  Assume $g : \sX \rightarrow \mathbb{R}$ such that $\text{E}_{\pi}
  |g|^{2+\delta} < \infty$ for some $\delta > 0$ and let $X$ be a
  Harris ergodic Markov chain with invariant distribution $\pi$.
  Further, suppose $X$ is uniformly ergodic.  If\vspace{-2mm}
\begin{enumerate}
\item $a_n \rightarrow \infty$ as $n \rightarrow \infty$,\vspace{-1mm}
\item $b_{n} \rightarrow \infty$ and $b_{n} / n \rightarrow 0$ as
  $n \rightarrow \infty$, \vspace{-1mm}
\item $b_{n}^{-1} n^{1-2\alpha} \log n \rightarrow 0$ as $n
  \rightarrow \infty$ where $\alpha \in (0, \delta/(24+12\delta))$
  and\vspace{-1mm}
\item there exists a constant $c \ge 1$ such that $\sum_{n} (b_{n} /
  n)^{c} < \infty$\vspace{-2mm}
\end{enumerate}
then as $n \rightarrow \infty$, $\hat{\sigma}_{BM}^{2} \rightarrow
\sigma_{g}^{2}$ w. p. 1.
\end{proposition}

In Appendix~\ref{app:sc_bm} we use part 2 of Lemma~\ref{lemma:sip} to
extend Proposition~\ref{prop:ue_sc_bm} to geometrically ergodic Markov
chains.

\begin{proposition} \label{prop:sc_bm}
  Assume $g : \sX \rightarrow \mathbb{R}$ such that $\text{E}_{\pi}
  |g|^{4+\delta} < \infty$ for some $\delta > 0$ and let $X$ be a
  Harris ergodic Markov chain with invariant distribution $\pi$.
  Further, suppose $X$ is geometrically ergodic.  If\vspace{-2mm}
\begin{enumerate}
\item $a_n \rightarrow \infty$ as $n \rightarrow \infty$,\vspace{-1mm}
\item $b_{n} \rightarrow \infty$ and $b_{n} / n \rightarrow 0$ as
  $n \rightarrow \infty$, \vspace{-1mm}
\item $b_{n}^{-1} n^{2\alpha} [ \log n ]^{3} \rightarrow 0$ as $n
  \rightarrow \infty$ where $\alpha = 1/(2 + \delta)$ and\vspace{-1mm}
\item there exists a constant $c \ge 1$ such that $\sum_{n} (b_{n} /
  n)^{c} < \infty$\vspace{-2mm}
\end{enumerate}
then as $n \rightarrow \infty$, $\hat{\sigma}_{BM}^{2} \rightarrow
\sigma_{g}^{2}$ w. p. 1.
\end{proposition}

\begin{remark}
There is no assumption of stationarity in
Propositions~\ref{prop:ue_sc_bm} or~\ref{prop:sc_bm}.  Hence burn-in
is not required to implement CBM.
\end{remark}

\begin{remark} \label{rem:bs}
  Consider using $b_{n} = \lfloor n^{\theta} \rfloor$ and
  $a_{n}=\lfloor n/b_{n} \rfloor$.  Proposition~\ref{prop:ue_sc_bm}
  requires that $1 > \theta > 1 - 2 \alpha > 1- \delta/(12+6\delta) > 5/6$
  but Proposition~\ref{prop:sc_bm} requires only $1 > \theta >
  (1+\delta/2)^{-1} >0$.
\end{remark}

Under the conditions of Propositions~\ref{prop:ue_sc_bm}
or~\ref{prop:sc_bm} an asymptotically valid fixed-width procedure for
estimating $E_{\pi} g$ results if we terminate the simulation the
first time
\[
t_{a_{n}-1} \frac{\hat{\sigma}_{BM}}{\sqrt{n}} + p(n)\le \epsilon
\]
where $t_{a_{n}-1}$ is the appropriate quantile from a student's $t$
distribution with $a_{n}-1$ degrees of freedom.

\subsection{Practical Implementation Issues}
\label{sec:practical}
Making practical use of the preceding theory requires (i) a moment
condition; (ii) establishing geometric ergodicity of the sampler at
hand; (iii) choosing $p(n)$; (iv) using RS requires \eqref{eq:mc} or
at least \eqref{eq:reg_pr}; and (v) CBM requires choosing $a_{n}$ and
$b_{n}$.

Since a moment condition is required even in the iid case we do not
view (i) as restrictive.  Consider (ii).  It is easy to construct
examples where the convergence rate is so slow that a Markov chain CLT
does not hold \citep{robe:1999} so the importance of establishing the
rate of convergence in \eqref{eq:tvbd} should not be underestimated.
On the other hand, the MCMC community has expended considerable effort
in trying to understand when certain Markov chains are geometrically
ergodic; see the references in Remark~\ref{rm:geo}.  In our view, this
is not the obstacle that it once was.

Regarding (iii), we know of no work on choosing an optimal $p(n)$.
Recall that the theory requires $p(n)=o(n^{-1/2})$.  In our examples
we use $p(n)=\epsilon I(n \le n^{*})$ where $n^{*} >0$ is fixed.
Since $n^{*}$ is typically chosen based on empirical experience with
the sampler at hand we might want a penalty for sample sizes greater
than $n^{*}$ so another reasonable choice might be $p(n)=\epsilon I(n
\le n^{*}) + C n^{-k}$ for some $k > 1/2$ and $C >0$.

The issue in (iv), i.e., calculating \eqref{eq:mc} or \eqref{eq:reg_pr}
is commonly viewed as overly burdensome.  However, in our experience,
this calculation need not be troublesome.  For example,
\citet{mykl:tier:yu:1995} give recipes for constructing \eqref{eq:mc}
and \eqref{eq:reg_pr} for Metropolis-Hastings independence and random
walk samplers; recall \eqref{eq:indep_rp}.  There is also some work on
establishing these conditions for very general models; see
\citet{hobe:jone:robe:2005}.  Finally, \citet{broc:kada:2005} and
\citet{geye:thom:1995} have shown that regenerations can be made to
occur naturally via simulated tempering.

Consider (v).  As we noted in Remark~\ref{rem:bs}, it is common to
choose the batch sizes according to $b_{n} = \lfloor n^{\theta}
\rfloor$ for some $\theta$.  \citet{song:schm:1995} and
\citet{chie:1988} have addressed the issue of what value of $\theta$
should be used from different theoretical points of view.  In
particular, \citet{chie:1988} showed that (under regularity
conditions) using $\theta = 1/2$ results in the batch means
approaching asymptotic normality at the fastest rate.
\citet{song:schm:1995} showed that (under different regularity
conditions) using $\theta=1/3$ minimizes the asymptotic mean-squared
error of $\hat{\sigma}^{2}_{BM}$. Note that Remark~\ref{rem:bs} shows
that $\theta=1/3$ requires a stronger moment condition than
$\theta=1/2$. We further address this issue in
Section~\ref{sec:examples}.

\subsection{Alternatives to BM and RS}  
\label{sec:alternative}
We chose to focus on BM and RS since in MCMC settings they seem to be
the most common methods for estimating the variance of the asymptotic
normal distribution.  However, there are other methods which may enjoy
strong consistency; e.g.  see \citet{dame:1991}, \citet{geye:1992},
\citet{numm:2002} and \citet{peli:shao:1995}.  In particular,
\citet{dame:1991} uses a strong invariance principle to obtain strong
consistency of certain spectral variance estimators under conditions
similar to those required in Proposition~\ref{prop:ue_sc_bm}.
Apparently, this can be extended to geometrically ergodic chains via
Lemma~\ref{lemma:sip} to obtain a result with regularity conditions
similar to Proposition~\ref{prop:sc_bm}.  However, we do not pursue
this further here.

\section{Examples}
\label{sec:examples}
In this section we investigate the finite sample performance of RS, BM
with 30 batches, and CBM with $b_{n}=\lfloor n^{1/3} \rfloor$ and
$b_{n}=\lfloor n^{1/2} \rfloor$ in four examples.  In particular, we
examine the coverage probabilities and half-widths of the resulting
intervals as well as the required simulation effort.  While each
example concerns a different statistical model and MCMC sampler there
are some commonalities.  In each case we perform many independent
replications of the given MCMC sampler.  The number of replications
ranges from 2000 to 9000 depending on the complexity of the example.
We used all methods on the \textit{same} output from each replication
of the MCMC sampler.  When the half-width of a 95\% interval with
$p(n) = \epsilon I(n \ge n^{*})$ (or $p(R) = \epsilon I(R \ge R^{*})$
for RS) is less than $\epsilon$ for a particular method, that
procedure was stopped and the chain length recorded.  Our choice of
$n^{*}$ is different for each example and was chosen based on our
empirical experience with the given Markov chain.  Other procedures
would continue until all of them were below the targeted half-width,
at which time a single replication was complete.  In order to estimate
the coverage probabilities we need true values of the quantities of
interest.  These are not analytically available in three of our
examples.  Our solution is to obtain precise estimates of the truth
through independent methods which are different for each example.  The
details are described below.  The results are reported in
Table~\ref{tab:summary}.

\subsection{Toy Example}
\label{sec:toy}
Consider estimating the mean of a $\text{Pareto}(\alpha, \beta)$
distribution, i.e., $\alpha \beta / (\beta-1)$, $\beta > 1$, using a
Metropolis-Hastings independence sampler with a $\text{Pareto}
(\alpha, \lambda)$ candidate.  Let $\pi$ be the target density and
$\nu$ be the proposal density.  Assume $\beta \ge \lambda$. Then for
$x \ge \alpha$
\[
\frac{\pi(x)}{\nu(x)} = \frac{\beta}{\lambda} \alpha^{\beta -
  \lambda}x^{\lambda-\beta} \le \frac{\beta}{\lambda} \; .
\]  
By Theorem 2.1 in \citet{meng:twee:1996} this sampler is uniformly
ergodic and
\[
\|P^{n}(x,\cdot) - \pi(\cdot)\| \le \left(1 -
  \frac{\lambda}{\beta}\right)^{n} \; .
\] 
In order to ensure the moment conditions required for
Proposition~\ref{prop:sc_bm} we set $\beta=10$ and $\lambda=9$ in
which case the right hand side is $10^{-n}$. Hence this sampler
converges extremely fast.  Implementation of RS was accomplished using
\eqref{eq:indep_rp} with $c=1.5$.

\subsubsection{Comparing convergence diagnostics with CBM}
\label{sec:compare}
As noted by a referee, one method for terminating the simulation is
via convergence diagnostics.  Consider the method of \citet{gewe:1992}
which is a diagnostic that seems close in spirit to the current work.
Geweke's diagnostic (GD) is based on a Markov chain CLT and hence does
not apply much more generally than CBM; the same can be said for many
other diagnostics.  GD uses a hypothesis test to ascertain when
$\bar{g}_{n}$ has stabilized.

In the remainder of this subsection we compare GD and CBM in terms of
mean-squared error (MSE) and chain length.  To this end we ran 9000
independent replications of the independence sampler with $\alpha=1$,
$\beta=10$ and $\lambda=9$.  We used CBM and GD on the output in the
following manner.  For each replication we set $n^{*}=45$ but the R
package {\tt boa} required a minimum of 120 iterations in order to
calculate GD.  After the minimum was achieved and the cutoff for a
particular method was attained we noted the chain length and the
current estimate of $E_{\pi} g$.  The cutoff for CBM was to set the
desired half-width to $\epsilon=.005$.  The result of using GD is a
p-value.  We chose four values (.05, .10, .2 and .4) for the threshold
in an attempt to tune the computation.  The results are reported in
Table~\ref{tab:pareto}.  As we previously noted, this sampler mixes
extremely well.  Thus it is not surprising that using GD results in a
small estimated MSE.  However, using CBM results in much smaller MSE
than GD.  The average chain lengths make it is clear that GD stops the
simulation much too soon. Moreover, changing the p-value threshold for
GD does not result in substantial improvements in estimation accuracy.

\begin{table}
\begin{tabular}{|c|c|c|c|}
   \hline
   Method  & Cutoff & Estimated MSE & Average Chain Length\\\hline
CBM ($b_{n}=\lfloor n^{1/3} \rfloor$) & $\epsilon=.005$ & $6.65 \times
   10^{-6} (9.9 \times 10^{-8})$ & 2428 (5) \\
CBM ($b_{n}=\lfloor n^{1/2} \rfloor$)& $\epsilon=.005$ & $7.34 \times
   10^{-6} (1.2 \times 10^{-8})$ & 2615 (3)     \\ \hline 
Geweke     & p-value=.4 & $1.17 \times 10^{-4} (2\times 10^{-6})$
   &202.6 (3.4)\\ 
Geweke     & p-value=.2 & $1.30 \times 10^{-4} (2\times 10^{-6})$
   &148.9 (1.6) \\ 
Geweke     & p-value=.1 & $1.34 \times 10^{-4} (2\times 10^{-6})$
   &133.4 (.9) \\ 
Geweke     & p-value=.05& $1.37 \times 10^{-4} (2\times 10^{-6})$
   &127.4 (.5)  \\ \hline 
\end{tabular}
\caption{\label{tab:pareto} Summary statistics for CBM versus GD for
   Example 4.1. Standard errors of estimates are in parentheses.} 
\end{table}

\subsection{A Hierarchical Model}
\label{sec:hm}
\citet{efro:morr:1975} present a data set that gives the raw batting
averages (based on 45 official at-bats) and a transformation
($\sqrt{45} \, \text{arcsin}(2x-1)$) for 18 Major League Baseball
players during the 1970 season.  \citet{rose:1996} considers the
following conditionally independent hierarchical model for the
transformed data.  Suppose for $i=1,\ldots,K$ that
\begin{eqnarray}
\label{eq:rose_model}
Y_{i} | \theta_{i} \sim  \mbox{N}(\theta_{i},1) & \hspace*{5mm} &
\theta_{i} | \mu, \lambda  \sim  \mbox{N}(\mu, \lambda)  \\
\lambda \sim \mbox{IG}(2, 2) & & f(\mu) \propto  1 \; .\nonumber
\end{eqnarray}
(Note that we say $W \sim \text{Gamma} (\alpha, \beta)$ if its density
is proportional to $w^{\alpha - 1} e^{-\beta w} I(w > 0)$ and if $X
\sim \text{Gamma}(b,c)$ then $X^{-1} \sim \text{IG}(b,c)$.)
\citet{rose:1996} introduces a Harris ergodic block Gibbs sampler that
has the posterior, $\pi(\theta,\mu,\lambda|y)$, characterized by the
hierarchy in \eqref{eq:rose_model} as its invariant distribution. This
Gibbs sampler completes a one-step transition $(\lambda', \mu',
\theta') \rightarrow (\lambda, \mu, \theta)$ by drawing from the
distributions of $\lambda | \theta'$ then $\mu | \theta', \lambda$ and
subsequently $\theta | \mu, \lambda$.  The full conditionals needed to
implement this sampler are given by
$$
\lambda | \theta, y \sim \mathrm{IG} \left( 2 + \frac{K-1}{2}, 
2 + \frac{\sum{(\theta_i - \bar{\theta})^{2}}}{2} \right), \hspace*{3mm}
\mu | \theta, \lambda, y \sim \mathrm{N} \left( \bar{\theta},
\frac{\lambda}{K} \right),
$$
$$
\theta_i | \lambda, \mu, y \stackrel{\mathrm{ind}}{\sim}~ \mathrm{N} 
\left( \frac{\lambda y_i +  \mu}{\lambda + 1}, \frac{\lambda }
{\lambda + 1} \right).
$$
Rosenthal proved geometric ergodicity of the associated Markov
chain.  However, MCMC is not required to sample from the posterior; in
Appendix~\ref{sec:rose} we develop an accept-reject sampler that
produces an iid sample from the posterior.   Also in
Appendix~\ref{sec:rose} we derive an expression for the probability of
regeneration \eqref{eq:reg_pr}. 

We focus on estimating the posterior mean of $\theta_{9}$, the
``true'' long-run (transformed) batting average of the Chicago Cubs'
Ron Santo. It is straightforward to check that the moment conditions
for CBM and RS are met.  Finally, we employed our accept-reject
sampling algorithm to generate $9 \times 10^{7}$ independent draws
from $\pi(\theta_{9} |y)$ which were then used to estimate the
posterior mean of $\theta_{9}$ which we assumed to be the truth.

\subsection{Calculating Exact Conditional P-Values}
\label{sec:p-value}

\citet[][p. 432]{agre:2002} reports data that correspond to pairs of
scorings of tumor ratings by two pathologists. A linear by linear
association model specifies that the log of the Poisson mean in cell
$i,j$ satisfies
$$
\log \mu_{ij} = \alpha + \beta_i + \gamma_j + \delta \, i j \; .
$$
A parameter free null distribution for testing goodness-of-fit is
obtained by conditioning on the sufficient statistics for the
parameters, i.e., the margins of the table and $\sum_{ij} n_{ij} \,
ij$, where the $n_{ij}$ are the observed cell counts.  The resulting
conditional distribution is a generalization of the hypergeometric
distribution. An exact p-value for goodness-of-fit versus a saturated
alternative can be calculated by summing the conditional probabilities
of all tables satisfying the margins and the additional constraint and
having deviance statistics larger than the observed. 

For the current data set there are over twelve billion tables that
satisfy the margin constraints but an exhaustive search revealed that
there are only roughly 34,000 tables that also satisfy the constraint
induced by $\sum_{ij} n_{ij} \, ij$.  We will denote this set of
permissible tables by $\Gamma$.  Now the desired p-value is given by
\begin{equation}
\label{eq:p-value}
\sum_{y \in \Gamma} I[d(y) \ge d(y_{obs})] \, \pi(y)
\end{equation}
where $d(\cdot)$ is the deviance function and $\pi$ denotes the
generalized hypergeometric.  Since we have enumerated $\Gamma$ we find
that the true exact p-value is .044 whereas the chi-squared
approximation yields a p-value of .368.  However, a different data set
with different values of the sufficient statistics will have a
different reference set which must be enumerated in order to find the
exact p-value.  This would be too computationally burdensome to
implement generally and hence it is common to resort to MCMC-based
approximations \citep[see e.g.][]{caff:boot:2001,diac:stur:1998}.

To estimate \eqref{eq:p-value} we will use the Metropolis-Hastings
algorithm developed in \citet{caff:boot:2001}.  This algorithm is also
employed by the R package \texttt{exactLoglinTest}.  The associated
Markov chain is Harris ergodic and its invariant distribution is the
appropriate generalized hypergeometric distribution.  Moreover, the
chain is uniformly ergodic and since we are estimating the expectation
of a bounded function the regularity conditions for both RS and CBM
are easily met.

Implementation of RS is straightforward.  As we mentioned earlier, in
finite state spaces regenerations occur whenever the chain returns to
any fixed state.  In order to choose the fixed state we ran the
algorithm for 1000 iterations and chose the state which had the
highest probability with respect to the stationary distribution.  The
same fixed state was used in each replication. 

\subsection{A Model-Based Spatial Statistics Application}
\label{sec:spatial}
Consider the Scottish lip cancer data set \citep{clay:kald:1987} which
consists of the number of cases of lip cancer registered in each of
the 56 (pre-reorganization) counties of Scotland, together with the
expected number of cases given the age-sex structure of the
population.  We assume a Poisson likelihood for areal (spatially
aggregated) data.  Specifically, for $i=1,....,N$ we assume that given
$\mu_i$ the disease counts $Y_{i}$ are conditionally independent and
\begin{equation}
Y_{i} | \mu_{i} \sim \text{Poisson}(E_{i}e^{\mu_i})
\end{equation}
where $E_{i}$ is the known `expected' number of disease events in the
$i$th region assuming constant risk and $\mu_i$ is the log-relative
risk of disease for the $i$th region.  Set $\phi = (\phi_{1}, \ldots,
\phi_{N})^{T}$.  Each $\mu_i$ is modeled as $\mu_i = \theta_i +
\phi_i$ where
\begin{equation*}
    \theta_i|\tau_h \sim  \; \text{N} (0,1/\tau_h), \hspace*{6mm} 
    \phi | \tau_c \sim  \; \text{CAR} (\tau_c) \; \propto \; 
    \tau_c^{N/2}\exp\left(-\frac{\tau_c}{2}\phi^TQ\phi\right) , \text{
    and }
\end{equation*}
$$Q_{ij}=\left\{
\begin{array}{l l}
 \phantom{-}n_i & \mbox{if }i=j\\
 \phantom{-}0 & \mbox{if }i\mbox{ is not adjacent to }j\\
-1 & \mbox{if }i\mbox{ is adjacent to }j\\
\end{array} \right.
$$
with $n_i$ the number of neighbors for the $i$th region. Each
$\theta_i$ captures the $i$th region's extra-Poisson variability due
to area-wide heterogeneity, while each $\phi_i$ captures the $i$th
region's excess variability attributable to regional clustering. The
priors on the precision parameters are $\tau_h \sim \text{Gamma} (1,
.01)$ and $\tau_c \sim \text{Gamma} (1, .02)$.  This is a challenging
model to consider since the random effects parameters
($\theta_i,\phi_i$) are not identified in the likelihood, and the
spatial prior used is improper.  Also, no closed form expressions are
available for the marginal distributions of the parameters, and the
posterior distribution has $2N+2$ dimensions (114 for the lip cancer
data).

\citet{hara:tier:2004} establish uniform ergodicity of a Harris
ergodic Metropolis-Hastings independence sampler with invariant
distribution $\pi(\theta, \phi, \tau_h, \tau_c | y)$ where
$\theta=(\theta_{1}, \ldots, \theta_{N})^{T}$ and a heavy-tailed
proposal.  In our implementation of RS we used the formula for the
probability of a regeneration given by \eqref{eq:indep_rp} with $\log
c= -342.72$.  Using the empirical supremum of the ratio of the
invariant density to the proposal density (based on several draws from
the proposal) guided the choice of $c$.

We focus on estimating the posterior expectation of $\phi_{7}$, the
log-relative risk of disease for County 7 attributable to spatial
clustering.   Finally, we used an independent run of length
$10^{7}$ to obtain an estimate which we treated as the `true value'.

\subsection{Summary}
\label{sec:summary}
Table~\ref{tab:summary} reveals that the estimates of the coverage
probabilities are all less than the desired .95.  However, examining
the standard errors shows that only BM is significantly less in all of
the examples and the estimated coverage probability for RS is
\textit{not} significantly different from .95 in 3 out of 4.  The
story for CBM is more complicated in that the coverage depends on the
choice of $b_{n}$.  Using $b_{n}=\lfloor n^{1/3} \rfloor$ gives the
best coverage for the examples in Sections~\ref{sec:toy}
and~\ref{sec:hm} while $b_{n}=\lfloor n^{1/2} \rfloor$ is superior for
those in Sections~\ref{sec:p-value} and~\ref{sec:spatial}.  The reason
for this is that the Markov chains in Sections~\ref{sec:toy}
and~\ref{sec:hm} mix exceptionally well and hence smaller batch sizes
can be tolerated.  However, the examples in Sections~\ref{sec:p-value}
and~\ref{sec:spatial} are realistic problems and hence the chains do
not mix as well so that larger batch sizes are required.  Thus we
would generally recommend using $b_{n}=\lfloor n^{1/2} \rfloor$.

The example in subsection~\ref{sec:p-value} deserves to be singled out
due to the low estimated coverage probabilities.  The goal in this
example was to estimate a fairly small probability, a situation in
which the Wald interval is known to have poor coverage even in iid
settings.

While RS and CBM appear comparable in terms of coverage probability RS
tends to result in slightly longer runs than CBM which in turn results
in longer runs than BM.  Moreover, RS and CBM are comparable in their
ability to produce intervals that meet the target half-width more
closely than BM.  Also, the intervals for RS are apparently more
stable than those of CBM and BM.  Finally, BM underestimates the Monte
Carlo standard error and therefore suggests stopping the chain too
early.

While RS has a slight theoretical advantage over CBM their finite
sample properties appear comparable.  Also, like RS, CBM avoids the
burn-in issue, which has been a long standing obstacle to MCMC
practitioners. In addition, CBM enjoys the advantage of being slightly
easier to implement.  Thus CBM clearly has a place in the tool kit of
MCMC practitioners.

\begin{appendix}
\section{Proof of Lemma~\ref{lemma:rsmoments}}
\label{app:rsmoments}

\subsection{Preliminary Results}
\label{app:rs_prelim}
Recall the split chain $X'$ and that $0 = \tau_0 < \tau_1 < \tau_2 <
\cdots$ denote the regeneration times; i.e., $\tau_{r+1} = \min \{i >
\tau_{r} : \delta_{i-1}=1 \}$.

\begin{lemma} \citep[Lemma 1]{hobe:jone:pres:rose:2002}
  \label{lemma:hjpr1} 
  Let $X$ be a Harris ergodic Markov chain and assume that
  \eqref{eq:mc} holds.  Then for any function $h : \sX^{\infty}
  \rightarrow \mathbb{R}$
\[
\text{E}_{\pi} | h(X_{0}, X_{1}, \ldots )| \ge c \text{E}_{Q} |
h(X_{0}, X_{1}, \ldots )| 
\]
where $c = \text{E}_{\pi} s$. 
\end{lemma}

\begin{lemma} \citep[Lemma 2]{hobe:jone:pres:rose:2002}
  \label{lemma:hjpr2} 
  Let $X$ be a Harris ergodic Markov chain and assume that
  \eqref{eq:mc} holds.  If $X$ is geometrically ergodic, then there
  exists a $\beta > 1$ such that $\text{E}_{\pi} \beta^{\tau_{1}} <
  \infty$.
\end{lemma}
\begin{corollary} \label{cor:hjpr3}
Assume the conditions of Lemma~\ref{lemma:hjpr2}. For any $a > 0$
\[
\sum_{i=0}^{\infty} \left[ \text{Pr}_{\pi} (\tau_{1} \ge i + 1)
  \right]^{a} \le \left(\text{E}_{\pi}
  \beta^{\tau_1}\right)^{a}\sum_{i=0}^{\infty} \beta^{-a(i+1)} <
  \infty \; .
\]
\end{corollary}

\subsection{Proof of Lemma~\ref{lemma:rsmoments}}
We prove only part 2 of the lemma as part 1 is similar.  Without loss
of generality we assume $0 < \delta < 1$.  By Lemma~\ref{lemma:hjpr1},
it is enough to verify that $\text{E}_{\pi}\tau_{1}^{p} < \infty$ and
$\text{E}_{\pi} S_{1}^{p+\delta} < \infty$.  Lemma~\ref{lemma:hjpr2}
shows that $\text{E}_{\pi}\tau_{1}^{p} < \infty$ for any $p>0$. Note
that
\begin{equation*}
\begin{split}
  & \left( \sum_{i=0}^{\tau_{1} - 1} g(X_{i})\right)^{p+\delta} \le
  \left( \sum_{i=0}^{\tau_{1} - 1} |g(X_{i})| \right)^{p+\delta} =
  \left( \sum_{i=0}^{\infty} I(0 \le i \le \tau_{1} - 1)
    |g(X_{i})| \right)^{p+\delta} \\
  & \le \sum_{i_{1}=0}^{\infty} \cdots \sum_{i_{p}=0}^{\infty}
  \sum_{i_{p+1}=0}^{\infty}\left[ \prod_{j=1}^{p} I(0 \le i_{j} \le
    \tau_{1} - 1) |g(X_{i_{j}})| \right]I(0 \le i_{p+1} \le \tau_{1} -
  1)|g(X_{i_{p+1}})|^{\delta}
\end{split}
\end{equation*}
and hence
\begin{equation*}
\begin{split}
   \text{E}_{\pi} S_{1}^{p+\delta} & \le \sum_{i_{1}=0}^{\infty} \cdots
  \sum_{i_{p}=0}^{\infty}
  \sum_{i_{p+1}=0}^{\infty}\text{E}_{\pi}\left( \left[
      \prod_{j=1}^{p+1} I(0 \le i_{j} \le \tau_{1} - 1)\right] \left[
      \prod_{j=1}^{p}|g(X_{i_{j}})| \right] |g(X_{i_{p+1}})|^{\delta}
  \right) \\
  & \le \sum_{i_{1}=0}^{\infty} \cdots \sum_{i_{p}=0}^{\infty}
  \sum_{i_{p+1}=0}^{\infty} \left[ \text{E}_{\pi} I(0 \le i_{1} \le
    \tau_{1} - 1) |g(X_{i_{1}})|^{2}
  \right]^{1/2} \times \\
  & \cdots \times \left[\text{E}_{\pi} I(0 \le i_{p} \le \tau_{1} - 1)
    |g(X_{i_{p}})|^{2^{p}} \right]^{1/2^{p}} \left[ \text{E}_{\pi} I(0
    \le i_{p+1} \le \tau_{1} - 1) |g(X_{i_{p+1}})|^{2^{p} \delta}
  \right]^{1/2^{p}}
\end{split}
\end{equation*}
where the second inequality follows with repeated application of
Cauchy-Schwartz.  Set $a_{j}=1+2^{j}/\delta$ and
$b_{j}=1+\delta/2^{j}$ for $j=1,2,\ldots,p$ and apply H{\"o}lder's
inequality to obtain 
\[
\text{E}_{\pi} I(0 \le i_{j} \le \tau_{1} - 1) |g(X_{i_{j}})|^{2^{j}}
\le \left[ \text{E}_{\pi} I(0 \le i_{j} \le \tau_{1} - 1)
\right]^{1/a_{j}} \left[ \text{E}_{\pi} |g(X_{i_{j}})|^{2^{j}+\delta}
\right]^{1/b_{j}} \; . 
\]
Note that 
\[
c_{j} := \left[ \left(\text{E}_{\pi}
    |g(X_{i_{j}})|^{2^{j}+\delta}\right)^{1/b_{j}} \right]^{1/2^{p}}<
\infty \; .
\]
Also, if $a_{p+1} = 1 + 2^{p}$ and $b_{p+1} = 1 + 1/2^{p}$ then
\[
\text{E}_{\pi} I(0 \le i_{p+1} \le \tau_{1} - 1)
|g(X_{i_{p+1}})|^{2^{p} \delta} \le \left[ \text{E}_{\pi} I(0 \le
  i_{p+1} \le \tau_{1} - 1) \right]^{\frac{1}{a_{p+1}}} \left[
  \text{E}_{\pi} |g(X_{i_{p+1}})|^{\delta(2^{p}+\delta)}
\right]^{\frac{1}{b_{p+1}}} \; .
\]
Notice that 
\[
c_{p+1} := \left[ \left(\text{E}_{\pi}
    |g(X_{i_{p+1}})|^{\delta(2^{p}+\delta)}\right)^{1/b_{j}}
\right]^{1/2^{p}} < \infty
\]
and set $c=\max\{c_{1},\ldots,c_{p+1}\}$.  Then an appeal to
Corollary~\ref{cor:hjpr3} yields 
\begin{equation*}
\begin{split}
  &\text{E}_{\pi} S_{1}^{p+\delta} \le c \left[ \prod_{j=1}^{p}
    \sum_{i_{j}=0}^{\infty} \{\Pr_{\pi}(\tau_{1} \ge i_{j} +
    1)\}^{1/(a_{j} 2^{j})} \right]\!\left[ \sum_{i_{p+1}=0}^{\infty}
    \{\Pr_{\pi}(\tau_{1} \ge i_{j} + 1)\}^{1/(a_{p+1} 2^{p})}\right]<
    \infty \; .
\end{split}
\end{equation*}

\section{Proof of Proposition~\ref{prop:sc_bm}}
\label{app:sc_bm}

\subsection{Preliminary Results}
\label{app:prelim} 
Recall that $B=\{B(t), t \ge 0\}$ denotes a standard Brownian motion.
Define
\begin{equation}
\label{eq:unbm}
\tilde{\sigma}_{*}^{2} = \frac{b_{n}}{a_{n} - 1} \sum_{j=0}^{a_{n} - 1}
\left( \bar{B}_{j}(b_{n}) - \bar{B}(n) \right)^{2}
\end{equation}
where $\bar{B}_{j}(b_n) = b_{n}^{-1} \left( B((j+1)b_{n}) - B(j b_{n})
\right)$ and $\bar{B}(n) = n^{-1} B(n)$.

\begin{lemma}\citep[][p. 508]{dame:1994}
  For all $\epsilon > 0$ and for almost all sample paths there exists
  $n_0 (\epsilon)$ such that for all $n \ge n_0$\vspace{-1mm}
\citep[][p. 508]{dame:1994}
\begin{equation}
\label{eq:bj_lil}
| \bar{B}_{j} (b_n) | \le \sqrt{2} (1 + \epsilon) b_{n}^{-1/2} [ \log
  (n / b_n) + \log \log n]^{1/2} \; .
\end{equation}
\end{lemma}
\vspace{-2mm}
\begin{lemma} \citep{csor:reve:1981}
For all $\epsilon > 0$ and for almost all sample
paths there exists $n_0 (\epsilon)$ such that for all $n \ge
n_0$\vspace{-1mm} 
\begin{equation}
\label{eq:b_lil}
|B(n)| < (1 + \epsilon) [2 n \log \log n]^{1/2} \; .
\end{equation}
\end{lemma}

\subsection{Proof of Proposition~\ref{prop:sc_bm}}

Proposition~\ref{prop:sc_bm} follows from Lemma~\ref{lemma:sip} and
the following two lemmas:
\begin{lemma} \citep[][Proposition 3.1]{dame:1994} \label{lemma:dbm}
Assume \vspace*{-2mm}
\begin{enumerate}
\item $b_{n} \rightarrow \infty$ and $n / b_{n} \rightarrow \infty$ as
  $n \rightarrow \infty$ and \vspace*{-2mm}
\item there exists a constant $c \ge 1$ such that $\sum_{n} (b_{n} /
  n)^{c} < \infty$ \vspace*{-2mm}
\end{enumerate}
then as $n \rightarrow \infty$, $\tilde{\sigma}_{*}^{2} \rightarrow
1$ a.s. 
\end{lemma}

\begin{lemma} \label{lemma:ours}
Assume that \eqref{eq:sip} holds with $\gamma(n)=n^{\alpha} \log n$
where $\alpha = 1/(2+\delta)$.  If \vspace*{-2mm}
\begin{enumerate}
\item $a_n \rightarrow \infty$ as $ n \rightarrow \infty$, \vspace*{-2mm}
\item $b_{n} \rightarrow \infty$ and $n / b_{n} \rightarrow \infty$ as
  $n \rightarrow \infty$ and \vspace*{-2mm}
\item $b_{n}^{-1} n^{2\alpha} [ \log n ]^{3} \rightarrow 0$ as $n
  \rightarrow \infty$ where $\alpha = 1/(2 + \delta)$  \vspace*{-2mm}
\end{enumerate} 
then as $n \rightarrow \infty$, $ \hat{\sigma}_{BM}^{2} - \sigma_{g}^{2}
\tilde{\sigma}_{*}^{2} \rightarrow 0$ a.s.
\end{lemma}

\noindent \textit{Proof of Lemma}~\ref{lemma:ours}.
Recall that $X=\{X_1 , X_2 , \ldots \}$ is a Harris ergodic Markov
chain.  Define the process $Y$ by $Y_{i} = g(X_{i})- \text{E}_{\pi}g$
for $i=1, 2, 3, \ldots$.  Then
\[
\hat{\sigma}_{BM}^{2} = \frac{b_n}{a_n -1} \sum_{j=0}^{a_n - 1} \left(
  \bar{Y}_{j} (b_n) - \bar{Y} (n) \right)^{2}
\]
where $\bar{Y}_{j} (b_n) = b_{n}^{-1} \sum_{i=1}^{b_n} Y_{jb_n + i}$
for $j=0,\ldots , a_{n} - 1$ and $\bar{Y} (n) = n^{-1} \sum_{i=1}^{n}
Y_{i}$. Since
\[
\bar{Y}_{j} (b_n) - \bar{Y}(n)  = \bar{Y}_{j} (b_n) - \bar{Y}(n) \pm
\sigma_{g} \bar{B}_{j} (b_n)
\, \pm \sigma_{g} \bar{B} (n)
\]
we have
\begin{equation*}
\begin{split}
  \left|\hat{\sigma}_{BM}^{2} - \sigma_{g}^{2} \tilde{\sigma}_{*}^{2}
  \right| & \le \frac{b_{n}}{a_{n} - 1} \sum_{j=0}^{a_{n} - 1} \left[
    (\bar{Y}_{j} (b_n) - \sigma_{g} \bar{B}_{j} (b_n) )^{2} + (\bar{Y} (n)
    - \sigma_{g}
    \bar{B}(n))^{2}  \right.\\
  & + |2(\bar{Y}_{j} (b_n) - \sigma_{g} \bar{B}_{j} (b_n)) (\bar{Y} (n) -
  \sigma_{g} \bar{B}(n))| + |2 \sigma_{g} (\bar{Y}_{j} (b_n) -
  \sigma_{g} \bar{B}_{j} (b_n))\bar{B}_{j} (b_n)| \\
  & +| 2\sigma_{g}(\bar{Y}_{j} (b_n) - \sigma_{g} \bar{B}_{j} (b_n))
  \bar{B}(n)| + | 2 \sigma_{g} (\bar{Y} (n) -
  \sigma_{g}\bar{B}(n))\bar{B}_{j}(b_n)| \\
  & \left. + |2 \sigma_{g} (\bar{Y} (n) - \sigma_{g} \bar{B}(n))\bar{B}(n)|
  \right] \; .
\end{split}
\end{equation*}
Now we will consider each term in the sum and show that it tends to
0. \vspace*{-2mm}
\begin{enumerate}
\item Our assumptions say that there exists a constant $C$ such that
  for all large $n$
\begin{equation}
\label{eq:asa1}
\left| \sum_{i=1}^{n} g(X_{i}) - n \text{E}_{\pi} g - \sigma_{g} B(n)
\right| < C n^{\alpha} \log n \hspace*{5mm} a.s.
\end{equation} 
Note that
\[
\bar{Y}_{j} (b_n) - \sigma_{g} \bar{B}_{j} (b_n) = \frac{1}{b_n} \left[
  \sum_{i=1}^{(j+1)b_n} Y_{i} - \sigma_{g} B((j+1)b_n)\right] -
\frac{1}{b_n} \left[ \sum_{i=1}^{jb_n} Y_{i} - \sigma_{g} B(jb_n)\right]
\]
and hence by \eqref{eq:asa1}
\begin{equation}
\label{eq:asa2}
\begin{split}
|\bar{Y}_{j} (b_n) - \sigma_{g} \bar{B}_{j} (b_n)|  & \le \frac{1}{b_n}
\left[ | \sum_{i=1}^{(j+1)b_n} Y_{i} - \sigma_{g} B((j+1)b_n)| + | 
\sum_{i=1}^{jb_n} Y_{i} - \sigma_{g} B(jb_n) | \right]  \\
& < \frac{2}{b_n} C n^{\alpha} \log n  
\end{split}
\end{equation}
Then
\[
\frac{b_{n}}{a_{n} - 1} \sum_{j=0}^{a_{n} - 1}(\bar{Y}_{j}
    (b_n) - \sigma_{g} \bar{B}_{j} (b_n) )^{2} < 4 C^{2} \frac{a_n}{a_n -
    1} b_{n}^{-1} n^{2\alpha} (\log n)^{2} \; \rightarrow 0
\]
as $n \rightarrow \infty$ by conditions 1 and 3.
 \vspace*{-2mm}
\item Apply \eqref{eq:asa1} to obtain
\begin{equation}
\label{eq:asa3}
| \bar{Y} (n) - \sigma_{g} \bar{B} (n) | = n^{-1} | \sum_{i=1}^{n}
  Y_{i} - \sigma_{g} B(n) | < C n^{\alpha - 1} \log n \; .
\end{equation}
Then
\[
\frac{b_{n}}{a_{n} - 1} \sum_{j=0}^{a_{n} - 1}(\bar{Y} (n) - \sigma_{g}
\bar{B}(n))^{2} < C^{2} \frac{a_n}{a_n - 1} \frac{b_n}{n} \frac{(\log
  n)^{2}}{n^{1-2\alpha}} \rightarrow 0
\]
as $n \rightarrow \infty$ by conditions 1 and 2 and since
$1 - 2\alpha > 0$. 
\vspace*{-2mm}
\item By \eqref{eq:asa2} and \eqref{eq:asa3}
\[
|2(\bar{Y}_{j} (b_n) - \sigma_{g} \bar{B}_{j} (b_n)) (\bar{Y} (n) - \sigma_{g}
\bar{B}(n))| < 4 C^{2} b_{n}^{-1} n^{2 \alpha - 1} (\log n)^{2} \; .
\]
Thus
\[
\frac{b_{n}}{a_{n} - 1} \sum_{j=0}^{a_{n} - 1}|2(\bar{Y}_{j} (b_n) -
\sigma_{g} \bar{B}_{j} (b_n)) (\bar{Y} (n) - \sigma_{g} \bar{B}(n))| < 4 C^{2}
\frac{a_n}{a_n - 1} \frac{(\log n)^{2}}{n^{1-2\alpha}} \rightarrow 0
\]
as $n \rightarrow \infty$ by condition 1 and since $1 - 2\alpha > 0$.
 \vspace*{-2mm}
\item Since $b_n \ge 2$, \eqref{eq:bj_lil} and \eqref{eq:asa2}
  together imply
\[
\begin{split}
  |(\bar{Y}_{j} (b_n) - \sigma_{g} \bar{B}_{j} (b_n))\bar{B}_{j} (b_n)| &
  < 2^{3/2} C (1 + \epsilon) b_{n}^{-1} \left[ b_{n}^{-1} n^{2 \alpha}
    (\log n)^{2} \log(n/b_{n}) \right. \\
  & \left. + b_{n}^{-1} n^{2 \alpha} (\log n)^{2} \log \log n
  \right]^{1/2}
\end{split}
\]
Hence
\[
\begin{split}
\frac{b_{n}}{a_{n} - 1} \sum_{j=0}^{a_{n} - 1}|2 \sigma_{g} (\bar{Y}_{j}
(b_n) - \sigma_{g} \bar{B}_{j} (b_n))\bar{B}_{j} (b_n)| & \le 8 \sigma_{g} C (1
+ \epsilon) \frac{a_n}{a_n - 1} \left[ b_{n}^{-1} n^{2 \alpha}
  (\log n)^{2} \log(n/b_{n}) \right. \\
& + \left. b_{n}^{-1} n^{2 \alpha} (\log n)^{2} \log \log n
\right]^{1/2} \; \rightarrow 0
\end{split}
\]
as $n \rightarrow \infty$ by conditions 1 and 3.
 \vspace*{-2mm}
\item  By \eqref{eq:asa2} and \eqref{eq:b_lil} $|(\bar{Y}_{j} (b_n) -
  \sigma_{g} \bar{B}_{j} (b_n)) \bar{B}(n)| < 4 C (1 + \epsilon)
  b_{n}^{-1} n^{-1/2 + \alpha} (\log n)(\log \log n)^{1/2}$ so that
\[
\frac{b_{n}}{a_{n} - 1} \sum_{j=0}^{a_{n} - 1}|2 \sigma_{g}(\bar{Y}_{j}
(b_n) - \sigma_{g} \bar{B}_{j} (b_n)) \bar{B}(n)| < 8 \sigma_{g} C (1 +
\epsilon) \frac{a_n}{a_n - 1} \frac{(\log n)(\log \log
  n)^{1/2}}{n^{1/2 - \alpha}} \; \rightarrow 0
\]
as $n \rightarrow \infty$ by condition 1 and since $1/2 - \alpha > 0$.
 \vspace*{-2mm}
\item Use \eqref{eq:bj_lil} and \eqref{eq:asa3} to get
\[
| (\bar{Y} (n) - \sigma_{g}\bar{B}(n))\bar{B}_{j}(b_n)| < \sqrt{2} C (1 +
\epsilon) \frac{n^{\alpha-1} \log n}{\sqrt{b_{n}}} \left[
  \log(n/b_{n}) + \log \log n \right]^{1/2}
\]
and hence using conditions 1, 2 and 3 shows that as $n \rightarrow
\infty$ 
\[
\begin{split}
  & \frac{b_{n}}{a_{n} - 1} \sum_{j=0}^{a_{n} - 1}|2 \sigma_{g} (\bar{Y}
  (n) - \sigma_{g}\bar{B}(n))\bar{B}_{j}(b_n)| < \\
  & 4 \sigma_{g} C (1 + \epsilon) \frac{a_n}{a_n - 1} \frac{b_n}{n}
   \left[ b_{n}^{-1} n^{2\alpha} ((\log n)^{2}\log(n/b_{n}) +
    (\log n)^{2} \log \log n )\right]^{1/2} \rightarrow 0
\end{split}
\]
 \vspace*{-3mm}
\item Now \eqref{eq:b_lil} and \eqref{eq:asa3} imply $| (\bar{Y} (n) -
  \sigma_{g} \bar{B}(n))\bar{B}(n)| < 2 C (1 + \epsilon) n^{-3/2 + \alpha}
  (\log n)^{3/2}$. Hence
\[
\frac{b_{n}}{a_{n} - 1} \sum_{j=0}^{a_{n} - 1}|2 \sigma_{g}(\bar{Y}
(n) - \sigma_{g} \bar{B}(n))\bar{B}(n)| < 4 \sigma_{g} C (1 +
\epsilon) \frac{a_n}{a_n - 1} \frac{b_n}{n} \frac{(\log
  n)^{3/2}}{n^{1/2 - \alpha}} \; \rightarrow 0
\]
as $n \rightarrow \infty$ by conditions 1 and 2 and since $1/2 -
\alpha > 0$. 
\end{enumerate} 

\section{Calculations for Example~\ref{sec:hm}}
\label{sec:rose}
We consider a slightly more general formulation of the model given in
\eqref{eq:rose_model}.  
Suppose for $i=1,\ldots,K$ that
\begin{eqnarray}
\label{eq:rose_model1}
Y_{i} | \theta_{i} \sim  \mbox{N}(\theta_{i},a) & \hspace*{5mm} &
\theta_{i} | \mu, \lambda  \sim  \mbox{N}(\mu, \lambda)  \\
\lambda \sim \mbox{IG}(b,c) & & f(\mu) \propto  1 \; .\nonumber
\end{eqnarray}
where $a, b, c$ are all known positive constants.

\subsection{Sampling from $\pi(\theta,\mu,\lambda|y)$}
\label{app:ss}
Let $\pi(\theta,\mu,\lambda|y)$ be the posterior distribution
corresponding to the hierarchy in \eqref{eq:rose_model1}.  Note
that $\theta$ is a vector containing all of the $\theta_{i}$ and that
$y$ is a vector containing all of the data.  Consider the
factorization
\begin{equation}
\label{eq:post}
\pi(\theta,\mu,\lambda|y) = \pi(\theta|\mu,\lambda,y)
\pi(\mu|\lambda,y) \pi(\lambda|y) .
\end{equation}
If it is possible to sequentially simulate from each of the densities
on the right-hand side of \eqref{eq:post} we can produce iid draws
from the posterior.  Now $\pi(\theta|\mu,\lambda,y)$ is the product of
independent univariate normal densities, i.e. $\theta_{i} | \mu ,
\lambda, y \sim \text{N}((\lambda y_{i} + a \mu)/(\lambda + a), \,
a\lambda / (\lambda + a))$.  Also, $\pi(\mu|\lambda,y)$ is a normal
distribution, i.e. $\mu | \lambda ,y \sim \text{N}(\bar{y}, (\lambda +
a)/K)$.  Next
\[
\pi(\lambda |y) \propto \frac{1}{\lambda^{b+1} (\lambda +
  a)^{(K-1)/2}} e^{-c/\lambda - s^{2}/2(\lambda +a)}
\]
where $\bar{y}=K^{-1} \sum_{i=1}^{K} y_{i}$ and $s^{2} =
\sum_{i=1}^{K} (y_{i} - \bar{y})^{2}$. An accept-reject sampler with
an $\text{IG}(b,c)$ candidate can be used to sample from $\pi(\lambda
|y)$ since if we let $g(\lambda)$ be the kernel of an $\text{IG}(b,c)$
  density
\[
\sup_{\lambda \ge 0} \frac{1}{g(\lambda) \lambda^{b+1} (\lambda +
 a)^{(K-1)/2}} e^{-c/\lambda - s^{2}/2(\lambda +a)} = \sup_{\lambda
 \ge 0} \, (\lambda + a)^{(1-K)/2} e^{-s^{2}/2(\lambda + a)} = M <
 \infty  
\]
It is easy to show that the only critical point is $\hat{\lambda}
=s^{2}/(K-1) - a$ which is where the maximum occurs if $\hat{\lambda}
> 0$.  But if $\hat{\lambda} \le 0$ then the maximum occurs at 0.

\subsection{Implementing regenerative simulation}
We begin by establishing the minorization condition \eqref{eq:mc} for
\pcite{rose:1996} block Gibbs sampler.  For the one-step transition
$(\lambda', \mu',\theta') \rightarrow (\lambda, \mu, \theta)$ the
Markov transition density, $p$, is given by $p(\lambda, \mu, \theta |
\lambda', \mu', \theta') = f(\lambda, \mu | \theta') f(\theta |
\lambda, \mu)$.  Note that $\mathsf{X} = \mathbb{R}^+ \times
\mathbb{R}^1 \times \mathbb{R}^K$.  Fix a point $(\tilde{\lambda},
\tilde{\mu}, \tilde{\theta}) \in \mathsf{X}$ and let $D \subseteq
\mathsf{X}$.  Then
\[
\begin{split}
  p(\lambda, \mu, \theta | \lambda', \mu', \theta')
  & = f(\lambda, \mu | \theta') f(\theta | \lambda, \mu) \\
  & \geq f(\lambda, \mu | \theta') f(\theta | \lambda, \mu)
  I_{\left\{(\lambda, \mu, \theta) \in D \right\}} \\
  & = \frac{f(\lambda, \mu | \theta')}{f(\lambda, \mu
    |\tilde{\theta})} f(\lambda, \mu | \tilde{\theta}) f(\theta |
  \lambda, \mu)
  I_{\left\{(\lambda, \mu, \theta) \in D \right\}} \\
  & \geq \left\{ \inf_{(\lambda, \mu, \theta) \in D} \frac{f(\lambda,
      \mu | \theta')}{f(\lambda, \mu |\tilde{\theta})} \right\}
  f(\lambda, \mu |\tilde{\theta}) f(\theta | \lambda, \mu)
  I_{\left\{(\lambda, \mu, \theta) \in D \right\}}
\end{split}
\]
and hence \eqref{eq:mc} will follow by setting
\[
\varepsilon  = \int_D f(\lambda, \mu|\tilde{\theta}) f(\theta | 
\lambda, \mu) ~d\lambda ~d\mu ~d\theta ,
\] 
\[
s(\lambda', \mu', \theta')  = \varepsilon \inf_{(\lambda, \mu, \theta) 
\in D} \frac{f(\lambda, \mu | \theta')}{f(\lambda, \mu
|\tilde{\theta})} \hspace*{5mm} \text{ and } \hspace*{5mm}
q(\lambda, \mu, \theta)  = \varepsilon^{-1} 
f(\lambda, \mu |\tilde{\theta})  f(\theta | \lambda, \mu) 
          I_{\left\{(\lambda, \mu, \theta) \in D \right\}}.
\]
Now using \eqref{eq:reg_pr} shows that when $(\lambda, \mu, \theta)
\in D$ the probability of regeneration is given by
\begin{equation}
\label{eq:prob}
\Pr(\delta=1 |\lambda', \mu', \theta', \lambda, \mu, \theta) = \left\{
  \inf_{(\lambda, \mu, \theta) \in D} \frac{f(\lambda, \mu |
    \theta')}{f(\lambda, \mu |\tilde{\theta})} \right\}
\frac{f(\lambda, \mu |\tilde{\theta})}{f(\lambda, \mu |\theta')}
\end{equation}
Thus we need to calculate the infimum and plug into \eqref{eq:prob}.
To this end let $0 < d_1 < d_2 < \infty$, $-\infty < d_3 < d_4 <
\infty$ and set $D = [d_1, d_2] \times [d_3, d_4] \times
\mathbb{R}^K$. Define $V(\theta, \mu) = \sum_{i=1}^K (\theta_i -
\mu)^{2}$ and note that
\[
\inf_{(\lambda, \mu, \theta) \in D} \frac{f(\lambda, \mu |
  \theta')}{f(\lambda, \mu |\tilde{\theta})} = \inf_{\lambda \in [d_1,
  d_2],~ \mu \in [d_3, d_4]} \exp \left\{\frac{V(\tilde{\theta}, \mu)
    - V(\theta', \mu)}{2 \lambda} \right\} = \exp \left\{
  \frac{V(\tilde{\theta}, \hat{\mu}) - V(\theta', \hat{\mu})}{2
    \hat{\lambda}} \right\}
\]
where $\hat{\mu} = d_{4} I(\bar{\theta'} \leq \bar{\tilde{\theta}}) +
d_{3} I(\bar{\theta'} > \bar{\tilde{\theta}})$ and $\hat{\lambda} =
d_{2}I(V(\theta', \hat{\mu}) \leq V(\tilde{\theta}, \hat{\mu})) +
d_{1} I(V(\theta', \hat{\mu}) > V(\tilde{\theta}, \hat{\mu}))$.  We
find the fixed point with a preliminary estimate of the mean of the
stationary distribution, and $D$ to be centered at that point.  Let
$(\tilde{\lambda}, \tilde{\mu}, \tilde{\theta})$ be the ergodic mean
for a preliminary Gibbs sampler run, and let $S_{\lambda}$ and
$S_{\mu}$ denote the usual sample standard deviations of $\lambda$ and
$\mu$ respectively.  After some trial and error we took $ d_1 = \max
\left\{.01,\tilde{\lambda} - .5 S_{\lambda}\right\}$, $d_2 =
\tilde{\lambda} + .5 S_{\lambda}$, $d_3 = \tilde{\mu} - S_{\mu}$ and
$d_4 = \tilde{\mu} + S_{\mu}$.
\end{appendix}

\bigskip
\bigskip
\noindent {\Large {\bf Acknowledgments}}
\bigskip

\noindent The authors are grateful to Ansu Chatterjee, Jeff Rosenthal
and Bill Sudderth for helpful conversations about this paper.

\begin{table}
\hspace{-10mm}\begin{tabular}{|c|c|c|c|c|c|c|}
   \hline
Example &  & & & Average half &Average Chain  & Coverage \\ 
Section &   Method & $b_{n}$ & $n^{*}$ / $R^{*}$  & width
   & Length & Probability \\  \hline 
4.1 & CBM  & $\lfloor n^{1/2} \rfloor$ & 45  & .0048 ($1.9\times
   10^{-6}$)& 2428 (5) & .923 (.003)\\ 
 9000 reps & CBM  & $\lfloor n^{1/3} \rfloor$ & 45 & .0049
   ($8.0\times 10^{-7}$)& 2615 (3) &.943 (.002)\\ 
$\epsilon=.005$ & BM & $\lfloor n/30 \rfloor$ & 45 & .0047
   $(2.4\times10^{-6})$ & 2342 (6) & .908 (.003) \\ 
&  RS & - & 30 & .0049 $(4.0\times 10^{-7})$ & 2653 (2) & .948
   (.002)\\ \hline
4.2 & CBM  & $\lfloor n^{1/2} \rfloor$ & 2000 & .0194
   ($7.2\times 10^{-6}$) & 5549 (13)     & .930 (.004) \\ 
5000 reps & CBM  & $\lfloor n^{1/3} \rfloor$ & 2000 & .0198 ($3.3
   \times 10^{-6}$) & 5778 (6) & .947 (.003) \\ 
$\epsilon=.02$ & BM & $\lfloor n/30 \rfloor$ & 2000 & .0191 ($1.1
   \times 10^{-5}$) & 5279 (18)& .915 (.004) \\ 
& RS   & - & 50 & .0198 ($2.3 \times 10^{-6}$) & 5818 (12)    &
   .945 (.003) \\ \hline 
4.3 & CBM  & $\lfloor n^{1/2} \rfloor$ & 4000 & .0049 ($1.6 \times
   10^{-6}$) & 56258 (405) & .920 (.006) \\   
2000 reps & CBM  & $\lfloor n^{1/3} \rfloor$&  4000 & .0049 ($1.8
   \times 10^{-6}$)& 46011 (499) & .869 (.008)\\ 
$\epsilon=.005$ & BM & $\lfloor n/30 \rfloor$ & 4000 & .0049 ($1.7
   \times 10^{-6}$) & 45768 (478) & .874 (.007) \\  
& RS   & - & 20 & .0049 ($4.3 \times 10^{-6}$) & 58265 (642)
   & .894 (.007) \\ \hline 
4.4 & CBM  & $\lfloor n^{1/2} \rfloor$ &  10000 & .00396
   ($8.0\times 10^{-7}$) & 168197 (270) & .934 (.005) \\  
2000 reps & CBM  & $\lfloor n^{1/3} \rfloor$ &  10000 & .00398
   ($4.0\times 10^{-7}$) & 137119 (125)& .900 (.006)\\  
$\epsilon=.004$ &  BM & $\lfloor n/30 \rfloor$ &  10000 & .00394
   ($1.2\times 10^{-6}$) & 132099 (809)   & .880 (.007)\\ 
&  RS  & - & 25 & .00398 ($2.0\times 10^{-7}$) & 179338
   (407)   & .942 (.005)\\ \hline 
\end{tabular}
\caption{\label{tab:summary} Summary statistics for BM, CBM and RS.
       Standard errors of estimates are in parentheses.}
\end{table}

\begin{thebibliography}{}

\bibitem[Agresti, 2002]{agre:2002}
Agresti, A. (2002).
\newblock {\em Categorical Data Analysis}.
\newblock Wiley, New York.

\bibitem[Billingsley, 1968]{bill:1968}
Billingsley, P. (1968).
\newblock {\em Convergence of Probability Measures}.
\newblock Wiley, New York.

\bibitem[Brockwell and Kadane, 2005]{broc:kada:2005}
Brockwell, A.~E. and Kadane, J.~B. (2005).
\newblock Identification of regeneration times in {MCMC} simulation, with
  application to adaptive schemes.
\newblock {\em Journal of Computational and Graphical Statistics}, 14:436--458.

\bibitem[Caffo and Booth, 2001]{caff:boot:2001}
Caffo, B.~S. and Booth, J.~G. (2001).
\newblock A {M}arkov chain {M}onte {C}arlo algorithm for approximating exact
  conditional probabilities.
\newblock {\em Journal of Compuatational and Graphical Statistics},
  10:730--745.

\bibitem[Chen, 1999]{chen:1999}
Chen, X. (1999).
\newblock Limit theorems for functionals of ergodic {M}arkov chains with
  general state space.
\newblock {\em Memoirs of the American Mathematical Society}, 139.

\bibitem[Chien, 1988]{chie:1988}
Chien, C.-H. (1988).
\newblock Small sample theory for steady state confidence intervals.
\newblock In Abrams, M., Haigh, P., and Comfort, J., editors, {\em Proceedings
  of the Winter Simulation Conference}, pages 408--413.

\bibitem[Chow and Robbins, 1965]{chow:robb:1965}
Chow, Y.~S. and Robbins, H. (1965).
\newblock On the asymptotic theory of fixed-width sequential confidence
  intervals for the mean.
\newblock {\em The Annals of Mathematical Statistics}, 36:457--462.

\bibitem[Christensen et~al., 2001]{chri:moll:waag:2001}
Christensen, O.~F., Moller, J., and Waagepetersen, R.~P. (2001).
\newblock Geometric ergodicity of {M}etropolis-{H}astings algorithms for
  conditional simulation in generalized linear mixed models.
\newblock {\em Methodology and Computing in Applied Probability}, 3:309--327.

\bibitem[Clayton and Kaldor, 1987]{clay:kald:1987}
Clayton, D. and Kaldor, J. (1987).
\newblock Empirical {B}ayes estimates of age-standardized relative risks for
  use in disease mapping.
\newblock {\em Biometrics}, 43:671--681.

\bibitem[Cowles and Carlin, 1996]{cowl:carl:1996}
Cowles, M.~K. and Carlin, B.~P. (1996).
\newblock Markov chain {M}onte {C}arlo convergence diagnostics: {A} comparative
  review.
\newblock {\em Journal of the American Statistical Association}, 91:883--904.

\bibitem[Cs{\'a}ki and Cs{\"o}rg\H{o}, 1995]{csak:csor:1995}
Cs{\'a}ki, E. and Cs{\"o}rg\H{o}, M. (1995).
\newblock On additive functionals of {M}arkov chains.
\newblock {\em Journal of Theoretical Probability}, 8:905--919.

\bibitem[Cs\"{o}rg\H{o} and R\'{e}v\'{e}sz, 1981]{csor:reve:1981}
Cs\"{o}rg\H{o}, M. and R\'{e}v\'{e}sz, P. (1981).
\newblock {\em Strong Approximations in Probability and Statistics}.
\newblock Academic Press.

\bibitem[Damerdji, 1991]{dame:1991}
Damerdji, H. (1991).
\newblock Strong consistency and other properties of the spectral variance
  estimator.
\newblock {\em Management Science}, 37:1424--1440.

\bibitem[Damerdji, 1994]{dame:1994}
Damerdji, H. (1994).
\newblock Strong consistency of the variance estimator in steady-state
  simulation output analysis.
\newblock {\em Mathematics of Operations Research}, 19:494--512.

\bibitem[Diaconis and Sturmfels, 1998]{diac:stur:1998}
Diaconis, P. and Sturmfels, B. (1998).
\newblock Algebraic algorithms for sampling from conditional distributions.
\newblock {\em The Annals of Statistics}, 26:363--397.

\bibitem[Douc and Soulier, 2004]{douc:fort:moul:soul:2004}
Douc, R., F. G. M.~E. and Soulier, P. (2004).
\newblock Practical drift conditions for subgeometric rates of convergence.
\newblock {\em The Annals of Applied Probability}, 14:1353--1377.

\bibitem[Doukhan et~al., 1994]{douk:mass:rio:1994}
Doukhan, P., Massart, P., and Rio, E. (1994).
\newblock The functional central limit theorem for strongly mixing processes.
\newblock {\em Annales de l'Institut Henri Poincare, Section B, Calcul des
  Probabilities et Statistique}, 30:63--82.

\bibitem[Efron and Morris, 1975]{efro:morr:1975}
Efron, B. and Morris, C. (1975).
\newblock Data analysis using {S}tein's estimator and its generalizations.
\newblock {\em Journal of the American Statistical Association}, 70:311--319.

\bibitem[Fishman, 1996]{fish:1996}
Fishman, G.~S. (1996).
\newblock {\em Monte Carlo: Concepts, Algorithms, and Applications}.
\newblock Springer, New York.

\bibitem[Fort and Moulines, 2000]{fort:moul:2000}
Fort, G. and Moulines, E. (2000).
\newblock V-subgeometric ergodicity for a {H}astings-{M}etropolis algorithm.
\newblock {\em Statistics and Probability Letters}, 49:401--410.

\bibitem[Fort and Moulines, 2003]{fort:moul:2003}
Fort, G. and Moulines, E. (2003).
\newblock Polynomial ergodicity of {M}arkov transition kernels.
\newblock {\em Stochastic Processes and their Applications}, 103:57--99.

\bibitem[Geweke, 1992]{gewe:1992}
Geweke, J. (1992).
\newblock Evaluating the accuracy of sampling-based approaches to the
  calculation of posterior moments (with discussion).
\newblock In Bernardo, J.~M., Berger, J.~O., Dawid, A.~P., and Smith, A. F.~M.,
  editors, {\em Bayesian Statistics 4. Proceedings of the Fourth Valencia
  International Meeting}, pages 169--188. Clarendon Press.

\bibitem[Geyer, 1992]{geye:1992}
Geyer, C.~J. (1992).
\newblock Practical {M}arkov chain {M}onte {C}arlo (with discussion).
\newblock {\em Statistical Science}, 7:473--511.

\bibitem[Geyer, 1999]{geye:1999}
Geyer, C.~J. (1999).
\newblock Likelihood inference for spatial point processes.
\newblock In Barndorff-Nielsen, O.~E., Kendall, W.~S., and van Lieshout, M.
  N.~M., editors, {\em Stochastic Geometry: Likelihood and Computation}, pages
  79--140. Chapman \& Hall/CRC, Boca {R}aton.

\bibitem[Geyer and Thompson, 1995]{geye:thom:1995}
Geyer, C.~J. and Thompson, E.~A. (1995).
\newblock Annealing {M}arkov chain {M}onte {C}arlo with applications to
  ancestral inference.
\newblock {\em Journal of the American Statistical Association}, 90:909--920.

\bibitem[Gilks et~al., 1998]{gilk:robe:sahu:1998}
Gilks, W.~R., Roberts, G.~O., and Sahu, S.~K. (1998).
\newblock Adaptive {M}arkov chain {M}onte {C}arlo through regeneration.
\newblock {\em Journal of the American Statistical Association}, 93:1045--1054.

\bibitem[Glynn and Iglehart, 1990]{glyn:igle:1990}
Glynn, P.~W. and Iglehart, D.~L. (1990).
\newblock Simulation output analysis using standardized time series.
\newblock {\em Mathematics of Operations Research}, 15:1--16.

\bibitem[Glynn and Whitt, 1991]{glyn:whit:1991}
Glynn, P.~W. and Whitt, W. (1991).
\newblock Estimating the asymptotic variance with batch means.
\newblock {\em Operations Research Letters}, 10:431--435.

\bibitem[Glynn and Whitt, 1992]{glyn:whit:1992}
Glynn, P.~W. and Whitt, W. (1992).
\newblock The asymptotic validity of sequential stopping rules for stochastic
  simulations.
\newblock {\em The Annals of Applied Probability}, 2:180--198.

\bibitem[Haran and Tierney, 2004]{hara:tier:2004}
Haran, M. and Tierney, L. (2004).
\newblock Perfect sampling for a {B}ayesian spatial model.
\newblock Technical report, Pennsylvania State University, Department of
  Statistics.

\bibitem[Hobert and Geyer, 1998]{hobe:geye:1998}
Hobert, J.~P. and Geyer, C.~J. (1998).
\newblock Geometric ergodicity of {G}ibbs and block {G}ibbs samplers for a
  hierarchical random effects model.
\newblock {\em Journal of Multivariate Analysis}, 67:414--430.

\bibitem[Hobert et~al., 2002]{hobe:jone:pres:rose:2002}
Hobert, J.~P., Jones, G.~L., Presnell, B., and Rosenthal, J.~S. (2002).
\newblock On the applicability of regenerative simulation in {M}arkov chain
  {M}onte {C}arlo.
\newblock {\em Biometrika}, 89:731--743.

\bibitem[Hobert et~al., 2005]{hobe:jone:robe:2005}
Hobert, J.~P., Jones, G.~L., and Robert, C.~P. (2005).
\newblock Using a {M}arkov chain to construct a tractable approximation of an
  intractable probability distribution.
\newblock {\em Scandinavian Journal of Statistics, {\rm to appear}}.

\bibitem[Ibragimov and Linnik, 1971]{ibra:linn:1971}
Ibragimov, I.~A. and Linnik, Y.~V. (1971).
\newblock {\em Independent and Stationary Sequences of Random Variables}.
\newblock Walters-Noordhoff, The Netherlands.

\bibitem[Jarner and Hansen, 2000]{jarn:hans:2000}
Jarner, S.~F. and Hansen, E. (2000).
\newblock Geometric ergodicity of {M}etropolis algorithms.
\newblock {\em Stochastic Processes and Their Applications}, 85:341--361.

\bibitem[Jarner and Roberts, 2002]{jarn:robe:2002}
Jarner, S.~F. and Roberts, G.~O. (2002).
\newblock Polynomial convergence rates of {M}arkov chains.
\newblock {\em Annals of Applied Probability}, 12:224--247.

\bibitem[Jones, 2004]{jone:2004}
Jones, G.~L. (2004).
\newblock On the {M}arkov chain central limit theorem.
\newblock {\em Probability Surveys}, 1:299--320.

\bibitem[Jones and Hobert, 2001]{jone:hobe:2001}
Jones, G.~L. and Hobert, J.~P. (2001).
\newblock Honest exploration of intractable probability distributions via
  {M}arkov chain {M}onte {C}arlo.
\newblock {\em Statistical Science}, 16:312--334.

\bibitem[Jones and Hobert, 2004]{jone:hobe:2004}
Jones, G.~L. and Hobert, J.~P. (2004).
\newblock Sufficient burn-in for {G}ibbs samplers for a hierarchical random
  effects model.
\newblock {\em The Annals of Statistics}, 32:784--817.

\bibitem[Lai, 2001]{lai:2001}
Lai, T.~L. (2001).
\newblock Sequential analysis: {S}ome classical problems and new challenges.
\newblock {\em Statistica Sinica}, 11:303--351.

\bibitem[Liu, 1997]{liu:1997}
Liu, W. (1997).
\newblock Improving the fully sequential sampling scheme of
  {A}nscombe-{C}how-{R}obbins.
\newblock {\em The Annals of Statistics}, 25:2164--2171.

\bibitem[Marchev and Hobert, 2004]{marc:hobe:2004}
Marchev, D. and Hobert, J.~P. (2004).
\newblock Geometric ergodicity of van {D}yk and {M}eng's algorithm for the
  multivariate {S}tudent's $t$ model.
\newblock {\em Journal of the American Statistical Association}, 99:228--238.

\bibitem[Mengersen and Tweedie, 1996]{meng:twee:1996}
Mengersen, K. and Tweedie, R.~L. (1996).
\newblock Rates of convergence of the {H}astings and {M}etropolis algorithms.
\newblock {\em The Annals of Statistics}, 24:101--121.

\bibitem[Meyn and Tweedie, 1993]{meyn:twee:1993}
Meyn, S.~P. and Tweedie, R.~L. (1993).
\newblock {\em Markov Chains and Stochastic Stability}.
\newblock Springer-Verlag, London.

\bibitem[Meyn and Tweedie, 1994]{meyn:twee:1994b}
Meyn, S.~P. and Tweedie, R.~L. (1994).
\newblock Computable bounds for geometric convergence rates of {M}arkov chains.
\newblock {\em The Annals of Applied Probability}, 4:981--1011.

\bibitem[Mira and Tierney, 2002]{mira:tier:2002}
Mira, A. and Tierney, L. (2002).
\newblock Efficiency and convergence properties of slice samplers.
\newblock {\em Scandinavian Journal of Statistics}, 29:1--12.

\bibitem[Mykland et~al., 1995]{mykl:tier:yu:1995}
Mykland, P., Tierney, L., and Yu, B. (1995).
\newblock Regeneration in {M}arkov chain samplers.
\newblock {\em Journal of the American Statistical Association}, 90:233--241.

\bibitem[Nadas, 1969]{nada:1969}
Nadas, A. (1969).
\newblock An extension of a theorem of {C}how and {R}obbins on sequential
  confidence intervals for the mean.
\newblock {\em The Annals of Mathematical Statistics}, 40:667--671.

\bibitem[Nummelin, 2002]{numm:2002}
Nummelin, E. (2002).
\newblock {MC}'s for {MCMC}'ists.
\newblock {\em International Statistical Review}, 70:215--240.

\bibitem[Peligrad and Shao, 1995]{peli:shao:1995}
Peligrad, M. and Shao, Q.-M. (1995).
\newblock Estimation of the variance of partial sums for $\rho$-mixing random
  variables.
\newblock {\em Journal of Multivariate Analysis}, 52:140--157.

\bibitem[Philipp and Stout, 1975]{phil:stou:1975}
Philipp, W. and Stout, W. (1975).
\newblock Almost sure invariance principles for partial sums of weakly
  dependent random variables.
\newblock {\em Memoirs of the American Mathematical Society}, 2:1--140.

\bibitem[Robert, 1995]{robe:1995}
Robert, C.~P. (1995).
\newblock Convergence control methods for {M}arkov chain {M}onte {C}arlo
  algorithms.
\newblock {\em Statistical Science}, 10:231--253.

\bibitem[Roberts, 1999]{robe:1999}
Roberts, G.~O. (1999).
\newblock A note on acceptance rate criteria for {CLT}s for
  {M}etropolis-{H}astings algorithms.
\newblock {\em Journal of Applied Probability}, 36:1210--1217.

\bibitem[Roberts and Polson, 1994]{robe:pols:1994}
Roberts, G.~O. and Polson, N.~G. (1994).
\newblock On the geometric convergence of the {G}ibbs sampler.
\newblock {\em Journal of the Royal Statistical Society, {\rm Series B}},
  56:377--384.

\bibitem[Roberts and Rosenthal, 1997]{robe:rose:1997c}
Roberts, G.~O. and Rosenthal, J.~S. (1997).
\newblock Geometric ergodicity and hybrid {M}arkov chains.
\newblock {\em Electronic Communications in Probability}, 2:13--25.

\bibitem[Roberts and Rosenthal, 1998]{robe:rose:1998b}
Roberts, G.~O. and Rosenthal, J.~S. (1998).
\newblock Markov chain {M}onte {C}arlo: {S}ome practical implications of
  theoretical results (with discussion).
\newblock {\em Canadian Journal of Statistics}, 26:5--31.

\bibitem[Roberts and Rosenthal, 1999]{robe:rose:1999a}
Roberts, G.~O. and Rosenthal, J.~S. (1999).
\newblock Convergence of slice sampler {M}arkov chains.
\newblock {\em Journal of the Royal Statistical Society, {\rm Series B}},
  61:643--660.

\bibitem[Roberts and Rosenthal, 2004]{robe:rose:2004}
Roberts, G.~O. and Rosenthal, J.~S. (2004).
\newblock General state space {M}arkov chains and {MCMC} algorithms.
\newblock {\em Probability Surveys}, 1:20--71.

\bibitem[Rosenthal, 1995]{rose:1995a}
Rosenthal, J.~S. (1995).
\newblock Minorization conditions and convergence rates for {M}arkov chain
  {M}onte {C}arlo.
\newblock {\em Journal of the American Statistical Association}, 90:558--566.

\bibitem[Rosenthal, 1996]{rose:1996}
Rosenthal, J.~S. (1996).
\newblock Analysis of the {G}ibbs sampler for a model related to
  {J}ames-{S}tein estimators.
\newblock {\em Statistics and Computing}, 6:269--275.

\bibitem[Song and Schmeiser, 1995]{song:schm:1995}
Song, W.~T. and Schmeiser, B.~W. (1995).
\newblock Optimal mean-squared-error batch sizes.
\newblock {\em Management Science}, 41:110--123.

\bibitem[Tierney, 1994]{tier:1994}
Tierney, L. (1994).
\newblock Markov chains for exploring posterior distributions (with
  discussion).
\newblock {\em The Annals of Statistics}, 22:1701--1762.

\end{thebibliography}
\end{document}